\documentstyle[12pt,bezier]{amsart}

\makeatletter
 
\def\diagram{\m@th\leftwidth=\z@ \rightwidth=\z@ \topheight=\z@
\botheight=\z@ \setbox\@picbox\hbox\bgroup}
 
\def\enddiagram{\egroup\wd\@picbox\rightwidth\unitlength
\ht\@picbox\topheight\unitlength \dp\@picbox\botheight\unitlength
\hskip\leftwidth\unitlength\box\@picbox}
 
\def\bfig{\begin{diagram}}
\def\efig{\end{diagram}}
\newcount\wideness \newcount\leftwidth \newcount\rightwidth
\newcount\highness \newcount\topheight \newcount\botheight
 
\def\ratchet#1#2{\ifnum#1<#2 \global #1=#2 \fi}
 
\def\putbox(#1,#2)#3{%
\horsize{\wideness}{#3} \divide\wideness by 2
{\advance\wideness by #1 \ratchet{\rightwidth}{\wideness}}
{\advance\wideness by -#1 \ratchet{\leftwidth}{\wideness}}
\vertsize{\highness}{#3} \divide\highness by 2
{\advance\highness by #2 \ratchet{\topheight}{\highness}}
{\advance\highness by -#2 \ratchet{\botheight}{\highness}}
\put(#1,#2){\makebox(0,0){$#3$}}}
 
\def\putlbox(#1,#2)#3{%
\horsize{\wideness}{#3}
{\advance\wideness by #1 \ratchet{\rightwidth}{\wideness}}
{\ratchet{\leftwidth}{-#1}}
\vertsize{\highness}{#3} \divide\highness by 2
{\advance\highness by #2 \ratchet{\topheight}{\highness}}
{\advance\highness by -#2 \ratchet{\botheight}{\highness}}
\put(#1,#2){\makebox(0,0)[l]{$#3$}}}
 
\def\putrbox(#1,#2)#3{%
\horsize{\wideness}{#3}
{\ratchet{\rightwidth}{#1}}
{\advance\wideness by -#1 \ratchet{\leftwidth}{\wideness}}
\vertsize{\highness}{#3} \divide\highness by 2
{\advance\highness by #2 \ratchet{\topheight}{\highness}}
{\advance\highness by -#2 \ratchet{\botheight}{\highness}}
\put(#1,#2){\makebox(0,0)[r]{$#3$}}}

\def\adjust[#1]{} 
 
\newcount \coefa
\newcount \coefb
\newcount \coefc
\newcount\tempcounta
\newcount\tempcountb
\newcount\tempcountc
\newcount\tempcountd
\newcount\xext
\newcount\yext
\newcount\xoff
\newcount\yoff
\newcount\gap%
\newcount\arrowtypea
\newcount\arrowtypeb
\newcount\arrowtypec
\newcount\arrowtyped
\newcount\arrowtypee
\newcount\height
\newcount\width
\newcount\xpos
\newcount\ypos
\newcount\run
\newcount\rise
\newcount\arrowlength
\newcount\halflength
\newcount\arrowtype
\newdimen\tempdimen
\newdimen\xlen
\newdimen\ylen
\newsavebox{\tempboxa}%
\newsavebox{\tempboxb}%
\newsavebox{\tempboxc}%
 
\newdimen\w@dth
 
\def\setw@dth#1#2{\setbox\z@\hbox{\m@th$#1$}\w@dth=\wd\z@
\setbox\@ne\hbox{\m@th$#2$}\ifnum\w@dth<\wd\@ne \w@dth=\wd\@ne \fi
\advance\w@dth by 1.2em}
 
\def\t@^#1_#2{\allowbreak\def\n@one{#1}\def\n@two{#2}\mathrel
{\setw@dth{#1}{#2}
\mathop{\hbox to \w@dth{\rightarrowfill}}\limits
\ifx\n@one\empty\else ^{\box\z@}\fi
\ifx\n@two\empty\else _{\box\@ne}\fi}}
\def\t@@^#1{\@ifnextchar_{\t@^{#1}}{\t@^{#1}_{}}}
\def\to{\@ifnextchar^{\t@@}{\t@@^{}}}
 
\def\t@left^#1_#2{\def\n@one{#1}\def\n@two{#2}\mathrel{\setw@dth{#1}{#2}
\mathop{\hbox to \w@dth{\leftarrowfill}}\limits
\ifx\n@one\empty\else ^{\box\z@}\fi
\ifx\n@two\empty\else _{\box\@ne}\fi}}
\def\t@@left^#1{\@ifnextchar_{\t@left^{#1}}{\t@left^{#1}_{}}}
\def\toleft{\@ifnextchar^{\t@@left}{\t@@left^{}}}
 
\def\two@^#1_#2{\allowbreak
\def\n@one{#1}\def\n@two{#2}\mathrel{\setw@dth{#1}{#2}
\mathop{\vcenter{\lineskip\z@\baselineskip\z@
                 \hbox to \w@dth{\rightarrowfill}%
                 \hbox to \w@dth{\rightarrowfill}}%
       }\limits
\ifx\n@one\empty\else ^{\box\z@}\fi
\ifx\n@two\empty\else _{\box\@ne}\fi}}
\def\tw@@^#1{\@ifnextchar _{\two@^{#1}}{\two@^{#1}_{}}}
\def\two{\@ifnextchar ^{\tw@@}{\tw@@^{}}}
 
\def\tofr@^#1_#2{\def\n@one{#1}\def\n@two{#2}\mathrel{\setw@dth{#1}{#2}
\mathop{\vcenter{\hbox to \w@dth{\rightarrowfill}\kern-1.7ex
                 \hbox to \w@dth{\leftarrowfill}}%
       }\limits
\ifx\n@one\empty\else ^{\box\z@}\fi
\ifx\n@two\empty\else _{\box\@ne}\fi}}
\def\t@fr@^#1{\@ifnextchar_ {\tofr@^{#1}}{\tofr@^{#1}_{}}}
\def\tofro{\@ifnextchar^ {\t@fr@}{\t@fr@^{}}}

\def\mon{\mathop{\m@th\hbox to
      14.6\P@{\lasyb\char'51\hskip-2.1\P@$\arrext$\hss
$\mathord\rightarrow$}}\limits} 
\def\leftmono{\mathrel{\m@th\hbox to
14.6\P@{$\mathord\leftarrow$\hss$\arrext$\hskip-2.1\P@\lasyb\char'50%
}}\limits} 
\mathchardef\arrext="0200

\def\settypes(#1,#2,#3){\arrowtypea#1 \arrowtypeb#2 \arrowtypec#3}
\def\settoheight#1#2{\setbox\@tempboxa\hbox{#2}#1\ht\@tempboxa\relax}%
\def\settodepth#1#2{\setbox\@tempboxa\hbox{#2}#1\dp\@tempboxa\relax}%
\def\settokens`#1`#2`#3`#4`{%
     \def\tokena{#1}\def\tokenb{#2}\def\tokenc{#3}\def\tokend{#4}}
\def\setsqparms[#1`#2`#3`#4;#5`#6]{%
\arrowtypea #1
\arrowtypeb #2
\arrowtypec #3
\arrowtyped #4
\width #5
\height #6
}
\def\setpos(#1,#2){\xpos=#1 \ypos#2}

\def\settriparms[#1`#2`#3;#4]{\settripairparms[#1`#2`#3`1`1;#4]}%
 
\def\settripairparms[#1`#2`#3`#4`#5;#6]{%
\arrowtypea #1
\arrowtypeb #2
\arrowtypec #3
\arrowtyped #4
\arrowtypee #5
\width #6
\height #6
}
 
\def\resetparms{\settripairparms[1`1`1`1`1;500]\width 500}
 
\resetparms
 
\def\mvector(#1,#2)#3{
\put(0,0){\vector(#1,#2){#3}}%
\put(0,0){\vector(#1,#2){26}}%
}
\def\evector(#1,#2)#3{{
\arrowlength #3
\put(0,0){\vector(#1,#2){\arrowlength}}%
\advance \arrowlength by-30
\put(0,0){\vector(#1,#2){\arrowlength}}%
}}
 
\def\horsize#1#2{%
\settowidth{\tempdimen}{$#2$}%
#1=\tempdimen
\divide #1 by\unitlength
}
 
\def\vertsize#1#2{%
\settoheight{\tempdimen}{$#2$}%
#1=\tempdimen
\settodepth{\tempdimen}{$#2$}%
\advance #1 by\tempdimen
\divide #1 by\unitlength
}
 
\def\putvector(#1,#2)(#3,#4)#5#6{{%
\ifnum3<\arrowtype
\putdashvector(#1,#2)(#3,#4)#5\arrowtype
\else
\ifnum\arrowtype<-3
\putdashvector(#1,#2)(#3,#4)#5\arrowtype
\else
\xpos=#1
\ypos=#2
\run=#3
\rise=#4
\arrowlength=#5
\ifnum \arrowtype<0
    \ifnum \run=0
        \advance \ypos by-\arrowlength
    \else
        \tempcounta \arrowlength
        \multiply \tempcounta by\rise
        \divide \tempcounta by\run
        \ifnum\run>0
            \advance \xpos by\arrowlength
            \advance \ypos by\tempcounta
        \else
            \advance \xpos by-\arrowlength
            \advance \ypos by-\tempcounta
        \fi
    \fi
    \multiply \arrowtype by-1
    \multiply \rise by-1
    \multiply \run by-1
\fi
\ifcase \arrowtype
\or \put(\xpos,\ypos){\vector(\run,\rise){\arrowlength}}%
\or \put(\xpos,\ypos){\mvector(\run,\rise)\arrowlength}%
\or \put(\xpos,\ypos){\evector(\run,\rise){\arrowlength}}%
\fi\fi\fi
}}
 
\def\putsplitvector(#1,#2)#3#4{
\xpos #1
\ypos #2
\arrowtype #4
\halflength #3
\arrowlength #3
\gap 140
\advance \halflength by-\gap
\divide \halflength by2
\ifnum\arrowtype>0
   \ifcase \arrowtype
   \or \put(\xpos,\ypos){\line(0,-1){\halflength}}%
       \advance\ypos by-\halflength
       \advance\ypos by-\gap
       \put(\xpos,\ypos){\vector(0,-1){\halflength}}%
   \or \put(\xpos,\ypos){\line(0,-1)\halflength}%
       \put(\xpos,\ypos){\vector(0,-1)3}%
       \advance\ypos by-\halflength
       \advance\ypos by-\gap
       \put(\xpos,\ypos){\vector(0,-1){\halflength}}%
   \or \put(\xpos,\ypos){\line(0,-1)\halflength}%
       \advance\ypos by-\halflength
       \advance\ypos by-\gap
       \put(\xpos,\ypos){\evector(0,-1){\halflength}}%
   \fi
\else \arrowtype=-\arrowtype
   \ifcase\arrowtype
   \or \advance \ypos by-\arrowlength
       \put(\xpos,\ypos){\line(0,1){\halflength}}%
       \advance\ypos by\halflength
       \advance\ypos by\gap
       \put(\xpos,\ypos){\vector(0,1){\halflength}}%
   \or \advance \ypos by-\arrowlength
       \put(\xpos,\ypos){\line(0,1)\halflength}%
       \put(\xpos,\ypos){\vector(0,1)3}%
       \advance\ypos by\halflength
       \advance\ypos by\gap
       \put(\xpos,\ypos){\vector(0,1){\halflength}}%
   \or \advance \ypos by-\arrowlength
       \put(\xpos,\ypos){\line(0,1)\halflength}%
       \advance\ypos by\halflength
       \advance\ypos by\gap
       \put(\xpos,\ypos){\evector(0,1){\halflength}}%
   \fi
\fi
}
 
\def\putmorphism(#1)(#2,#3)[#4`#5`#6]#7#8#9{{%
\run #2
\rise #3
\ifnum\rise=0
  \puthmorphism(#1)[#4`#5`#6]{#7}{#8}#9%
\else\ifnum\run=0
  \putvmorphism(#1)[#4`#5`#6]{#7}{#8}#9%
\else
\setpos(#1)%
\arrowlength #7
\arrowtype #8
\ifnum\run=0
\else\ifnum\rise=0
\else
\ifnum\run>0
    \coefa=1
\else
   \coefa=-1
\fi
\ifnum\arrowtype>0
   \coefb=0
   \coefc=-1
\else
   \coefb=\coefa
   \coefc=1
   \arrowtype=-\arrowtype
\fi
\width=2
\multiply \width by\run
\divide \width by\rise
\ifnum \width<0  \width=-\width\fi
\advance\width by60
\if l#9 \width=-\width\fi
\putbox(\xpos,\ypos){#4}
{\multiply \coefa by\arrowlength
\advance\xpos by\coefa
\multiply \coefa by\rise
\divide \coefa by\run
\advance \ypos by\coefa
\putbox(\xpos,\ypos){#5} }%
{\multiply \coefa by\arrowlength
\divide \coefa by2
\advance \xpos by\coefa
\advance \xpos by\width
\multiply \coefa by\rise
\divide \coefa by\run
\advance \ypos by\coefa
\if l#9%
   \putrbox(\xpos,\ypos){#6}%
\else\if r#9%
   \putlbox(\xpos,\ypos){#6}%
\fi\fi }%
{\multiply \rise by-\coefc
\multiply \run by-\coefc
\multiply \coefb by\arrowlength
\advance \xpos by\coefb
\multiply \coefb by\rise
\divide \coefb by\run
\advance \ypos by\coefb
\multiply \coefc by70
\advance \ypos by\coefc
\multiply \coefc by\run
\divide \coefc by\rise
\advance \xpos by\coefc
\multiply \coefa by140
\multiply \coefa by\run
\divide \coefa by\rise
\advance \arrowlength by\coefa
\ifcase\arrowtype
\or \put(\xpos,\ypos){\vector(\run,\rise){\arrowlength}}%
\or \put(\xpos,\ypos){\mvector(\run,\rise){\arrowlength}}%
\or \put(\xpos,\ypos){\evector(\run,\rise){\arrowlength}}%
\fi}\fi\fi\fi\fi}}

\newcount\numbdashes \newcount\lengthdash \newcount\increment
 
\def\howmanydashes{
\numbdashes=\arrowlength \lengthdash=40
\divide\numbdashes by \lengthdash
\lengthdash=\arrowlength
\divide\lengthdash by \numbdashes
\increment=\lengthdash
\multiply\lengthdash by 3
\divide\lengthdash by 5
}
 
\def\putdashvector(#1)(#2,#3)#4#5{%
\ifnum#3=0 \putdashhvector(#1){#4}#5
\else
\ifnum#2=0
\putdashvvector(#1){#4}#5\fi\fi}
 
\def\putdashhvector(#1,#2)#3#4{{%
\arrowlength=#3 \howmanydashes
\multiput(#1,#2)(\increment,0){\numbdashes}%
{\vrule height .4pt width \lengthdash\unitlength}
\arrowtype=#4 \xpos=#1
\ifnum\arrowtype<0 \advance\arrowtype by 7 \fi
\ifcase\arrowtype
\or \advance\xpos by 10
    \put(\xpos,#2){\vector(-1,0){\lengthdash}}
    \advance\xpos by 40
    \put(\xpos,#2){\vector(-1,0){\lengthdash}}
\or \advance \xpos by 10
    \put(\xpos,#2){\vector(-1,0){\lengthdash}}
    \advance\xpos by  \arrowlength
    \advance\xpos by  -50
    \put(\xpos,#2){\vector(-1,0){\lengthdash}}
\or \advance\xpos by 10
    \put(\xpos,#2){\vector(-1,0){\lengthdash}}
\or \advance\xpos by \arrowlength
    \advance\xpos by -\lengthdash
    \put(\xpos,#2){\vector(1,0){\lengthdash}}
\or {\advance\xpos by 10
    \put(\xpos,#2){\vector(1,0){\lengthdash}}}
    \advance\xpos by \arrowlength
    \advance\xpos by -\lengthdash
    \put(\xpos,#2){\vector(1,0){\lengthdash}}
\or \advance\xpos by \arrowlength
    \advance\xpos by -\lengthdash
    \put(\xpos,#2){\vector(1,0){\lengthdash}}
    \advance\xpos by -40
    \put(\xpos,#2){\vector(1,0){\lengthdash}}
   \fi
}}
 
\def\putdashvvector(#1,#2)#3#4{{%
\arrowlength=#3 \howmanydashes
\ypos=#2 \advance\ypos by -\arrowlength
\multiput(#1,#2)(0,\increment){\numbdashes}%
    {\vrule width .4pt height \lengthdash\unitlength}
\arrowtype=#4 \ypos=#2
\ifnum\arrowtype<0 \advance\arrowtype by 7 \fi
\ifcase\arrowtype
\or \advance\ypos by \arrowlength \advance\ypos by -40
    \put(#1,\ypos){\vector(0,1){\lengthdash}}
    \advance\ypos by -40
    \put(#1,\ypos){\vector(0,1){\lengthdash}}
\or \advance\ypos by 10
    \put(#1,\ypos){\vector(0,1){\lengthdash}}
    \advance\ypos by \arrowlength \advance\ypos by -40
    \put(#1,\ypos){\vector(0,1){\lengthdash}}
\or \advance\ypos by \arrowlength \advance\ypos by -40
    \put(#1,\ypos){\vector(0,1){\lengthdash}}
\or \advance\ypos by 10
    \put(#1,\ypos){\vector(0,-1){\lengthdash}}
\or \advance\ypos by 10
    \put(#1,\ypos){\vector(0,-1){\lengthdash}}
    \advance\ypos by \arrowlength \advance\ypos by -40
    \put(#1,\ypos){\vector(0,-1){\lengthdash}}
\or \advance\ypos by 10
    \put(#1,\ypos){\vector(0,-1){\lengthdash}}
    \advance\ypos by 40
    \put(#1,\ypos){\vector(0,-1){\lengthdash}}
\fi
}}
 
\def\puthmorphism(#1,#2)[#3`#4`#5]#6#7#8{{%
\xpos #1
\ypos #2
\width #6
\arrowlength #6
\arrowtype=#7
\putbox(\xpos,\ypos){#3\vphantom{#4}}%
{\advance \xpos by\arrowlength
\putbox(\xpos,\ypos){\vphantom{#3}#4}}%
\horsize{\tempcounta}{#3}%
\horsize{\tempcountb}{#4}%
\divide \tempcounta by2
\divide \tempcountb by2
\advance \tempcounta by30
\advance \tempcountb by30
\advance \xpos by\tempcounta
\advance \arrowlength by-\tempcounta
\advance \arrowlength by-\tempcountb
\putvector(\xpos,\ypos)(1,0)\arrowlength\arrowtype
\divide \arrowlength by2
\advance \xpos by\arrowlength
\vertsize{\tempcounta}{#5}%
\divide\tempcounta by2
\advance \tempcounta by20
\if a#8 %
   \advance \ypos by\tempcounta
   \putbox(\xpos,\ypos){#5}%
\else
   \advance \ypos by-\tempcounta
   \putbox(\xpos,\ypos){#5}%
\fi}}
 
\def\putvmorphism(#1,#2)[#3`#4`#5]#6#7#8{{%
\xpos #1
\ypos #2
\arrowlength #6
\arrowtype #7
\settowidth{\xlen}{$#5$}%
\putbox(\xpos,\ypos){#3}%
{\advance \ypos by-\arrowlength
\putbox(\xpos,\ypos){#4}}%
{\advance\arrowlength by-140
\advance \ypos by-70
\ifdim\xlen>0pt
   \if m#8%
      \putsplitvector(\xpos,\ypos)\arrowlength\arrowtype
   \else
   \putvector(\xpos,\ypos)(0,-1)\arrowlength\arrowtype
   \fi
\else
   \putvector(\xpos,\ypos)(0,-1)\arrowlength\arrowtype
\fi}%
\ifdim\xlen>0pt
   \divide \arrowlength by2
   \advance\ypos by-\arrowlength
   \if l#8%
      \advance \xpos by-40
      \putrbox(\xpos,\ypos){#5}%
   \else\if r#8%
      \advance \xpos by40
      \putlbox(\xpos,\ypos){#5}%
   \else
      \putbox(\xpos,\ypos){#5}%
   \fi\fi
\fi
}}
 
\def\putsquarep<#1>(#2)[#3;#4`#5`#6`#7]{{%
\setsqparms[#1]%
\setpos(#2)%
\settokens`#3`%
\puthmorphism(\xpos,\ypos)[\tokenc`\tokend`{#7}]{\width}{\arrowtyped}b%
\advance\ypos by \height
\puthmorphism(\xpos,\ypos)[\tokena`\tokenb`{#4}]{\width}{\arrowtypea}a%
\putvmorphism(\xpos,\ypos)[``{#5}]{\height}{\arrowtypeb}l%
\advance\xpos by \width
\putvmorphism(\xpos,\ypos)[``{#6}]{\height}{\arrowtypec}r%
}}
 
\def\putsquare{\@ifnextchar <{\putsquarep}{\putsquarep%
   <\arrowtypea`\arrowtypeb`\arrowtypec`\arrowtyped;\width`\height>}}
\def\square{\@ifnextchar< {\squarep}{\squarep
   <\arrowtypea`\arrowtypeb`\arrowtypec`\arrowtyped;\width`\height>}}
                                                   
\def\squarep<#1>[#2`#3`#4`#5;#6`#7`#8`#9]{{
\setsqparms[#1]
\diagram
\putsquarep<\arrowtypea`\arrowtypeb`\arrowtypec`
\arrowtyped;\width`\height>
(0,0)[#2`#3`#4`{#5};#6`#7`#8`{#9}]
\enddiagram
}}                                                 
                                                   
\def\putptrianglep<#1>(#2,#3)[#4`#5`#6;#7`#8`#9]{{%
\settriparms[#1]%
\xpos=#2 \ypos=#3
\advance\ypos by \height
\puthmorphism(\xpos,\ypos)[#4`#5`{#7}]{\height}{\arrowtypea}a%
\putvmorphism(\xpos,\ypos)[`#6`{#8}]{\height}{\arrowtypeb}l%
\advance\xpos by\height
\putmorphism(\xpos,\ypos)(-1,-1)[``{#9}]{\height}{\arrowtypec}r%
}}
 
\def\putptriangle{\@ifnextchar <{\putptrianglep}{\putptrianglep
   <\arrowtypea`\arrowtypeb`\arrowtypec;\height>}}
\def\ptriangle{\@ifnextchar <{\ptrianglep}{\ptrianglep
   <\arrowtypea`\arrowtypeb`\arrowtypec;\height>}}
                                              
\def\ptrianglep<#1>[#2`#3`#4;#5`#6`#7]{{
\settriparms[#1]
\diagram
\putptrianglep<\arrowtypea`\arrowtypeb`
\arrowtypec;\height>
(0,0)[#2`#3`#4;#5`#6`{#7}]
\enddiagram
}}                                            
 
\def\putqtrianglep<#1>(#2,#3)[#4`#5`#6;#7`#8`#9]{{%
\settriparms[#1]%
\xpos=#2 \ypos=#3
\advance\ypos by\height
\puthmorphism(\xpos,\ypos)[#4`#5`{#7}]{\height}{\arrowtypea}a%
\putmorphism(\xpos,\ypos)(1,-1)[``{#8}]{\height}{\arrowtypeb}l%
\advance\xpos by\height
\putvmorphism(\xpos,\ypos)[`#6`{#9}]{\height}{\arrowtypec}r%
}}
 
\def\putqtriangle{\@ifnextchar <{\putqtrianglep}{\putqtrianglep
   <\arrowtypea`\arrowtypeb`\arrowtypec;\height>}}
\def\qtriangle{\@ifnextchar <{\qtrianglep}{\qtrianglep
   <\arrowtypea`\arrowtypeb`\arrowtypec;\height>}}
                                              
\def\qtrianglep<#1>[#2`#3`#4;#5`#6`#7]{{
\settriparms[#1]
\width=\height                                
\diagram
\putqtrianglep<\arrowtypea`\arrowtypeb`
\arrowtypec;\height>
(0,0)[#2`#3`#4;#5`#6`{#7}]
\enddiagram
}}
 
\def\putdtrianglep<#1>(#2,#3)[#4`#5`#6;#7`#8`#9]{{%
\settriparms[#1]%
\xpos=#2 \ypos=#3
\puthmorphism(\xpos,\ypos)[#5`#6`{#9}]{\height}{\arrowtypec}b%
\advance\xpos by \height \advance\ypos by\height
\putmorphism(\xpos,\ypos)(-1,-1)[``{#7}]{\height}{\arrowtypea}l%
\putvmorphism(\xpos,\ypos)[#4``{#8}]{\height}{\arrowtypeb}r%
}}
 
\def\putdtriangle{\@ifnextchar <{\putdtrianglep}{\putdtrianglep
   <\arrowtypea`\arrowtypeb`\arrowtypec;\height>}}
\def\dtriangle{\@ifnextchar <{\dtrianglep}{\dtrianglep
   <\arrowtypea`\arrowtypeb`\arrowtypec;\height>}}
                                              
\def\dtrianglep<#1>[#2`#3`#4;#5`#6`#7]{{
\settriparms[#1]
\width=\height                                
\diagram
\putdtrianglep<\arrowtypea`\arrowtypeb`
\arrowtypec;\height>
(0,0)[#2`#3`#4;#5`#6`{#7}]
\enddiagram
}}
 
\def\putbtrianglep<#1>(#2,#3)[#4`#5`#6;#7`#8`#9]{{%
\settriparms[#1]%
\xpos=#2 \ypos=#3
\puthmorphism(\xpos,\ypos)[#5`#6`{#9}]{\height}{\arrowtypec}b%
\advance\ypos by\height
\putmorphism(\xpos,\ypos)(1,-1)[``{#8}]{\height}{\arrowtypeb}r%
\putvmorphism(\xpos,\ypos)[#4``{#7}]{\height}{\arrowtypea}l%
}}
 
\def\putbtriangle{\@ifnextchar <{\putbtrianglep}{\putbtrianglep
   <\arrowtypea`\arrowtypeb`\arrowtypec;\height>}}
\def\btriangle{\@ifnextchar <{\btrianglep}{\btrianglep
   <\arrowtypea`\arrowtypeb`\arrowtypec;\height>}}
                                             
\def\btrianglep<#1>[#2`#3`#4;#5`#6`#7]{{
\settriparms[#1]
\width=\height                               
\diagram
\putbtrianglep<\arrowtypea`\arrowtypeb`
\arrowtypec;\height>
(0,0)[#2`#3`#4;#5`#6`{#7}]
\enddiagram
}}
 
\def\putAtrianglep<#1>(#2,#3)[#4`#5`#6;#7`#8`#9]{{%
\settriparms[#1]%
\xpos=#2 \ypos=#3
{\multiply \height by2
\puthmorphism(\xpos,\ypos)[#5`#6`{#9}]{\height}{\arrowtypec}b}%
\advance\xpos by\height \advance\ypos by\height
\putmorphism(\xpos,\ypos)(-1,-1)[#4``{#7}]{\height}{\arrowtypea}l%
\putmorphism(\xpos,\ypos)(1,-1)[``{#8}]{\height}{\arrowtypeb}r%
}}
 
\def\putAtriangle{\@ifnextchar <{\putAtrianglep}{\putAtrianglep
   <\arrowtypea`\arrowtypeb`\arrowtypec;\height>}}
\def\Atriangle{\@ifnextchar <{\Atrianglep}{\Atrianglep
   <\arrowtypea`\arrowtypeb`\arrowtypec;\height>}}
                                                   
\def\Atrianglep<#1>[#2`#3`#4;#5`#6`#7]{{
\settriparms[#1]
\width=\height                                     
\diagram
\putAtrianglep<\arrowtypea`\arrowtypeb`
\arrowtypec;\height>
(0,0)[#2`#3`#4;#5`#6`{#7}]
\enddiagram
}}
 
\def\putAtrianglepairp<#1>(#2)[#3;#4`#5`#6`#7`#8]{{%
\settripairparms[#1]%
\setpos(#2)%
\settokens`#3`%
\puthmorphism(\xpos,\ypos)[\tokenb`\tokenc`{#7}]{\height}{\arrowtyped}b%
\advance\xpos by\height
\puthmorphism(\xpos,\ypos)[\phantom{\tokenc}`\tokend`{#8}]%
{\height}{\arrowtypee}b%
\advance\ypos by\height
\putmorphism(\xpos,\ypos)(-1,-1)[\tokena``{#4}]{\height}{\arrowtypea}l%
\putvmorphism(\xpos,\ypos)[``{#5}]{\height}{\arrowtypeb}m%
\putmorphism(\xpos,\ypos)(1,-1)[``{#6}]{\height}{\arrowtypec}r%
}}
 
\def\putAtrianglepair{\@ifnextchar <{\putAtrianglepairp}{\putAtrianglepairp%
   <\arrowtypea`\arrowtypeb`\arrowtypec`\arrowtyped`\arrowtypee;\height>}}
\def\Atrianglepair{\@ifnextchar <{\Atrianglepairp}{\Atrianglepairp%
   <\arrowtypea`\arrowtypeb`\arrowtypec`\arrowtyped`\arrowtypee;\height>}}
 
\def\Atrianglepairp<#1>[#2;#3`#4`#5`#6`#7]{{
\settripairparms[#1]
\settokens`#2`
\width=\height                                
\diagram
\putAtrianglepairp                            
<\arrowtypea`\arrowtypeb`\arrowtypec`
\arrowtyped`\arrowtypee;\height>
(0,0)[{#2};#3`#4`#5`#6`{#7}]
\enddiagram
}}
 
\def\putVtrianglep<#1>(#2,#3)[#4`#5`#6;#7`#8`#9]{{%
\settriparms[#1]%
\xpos=#2 \ypos=#3
\advance\ypos by\height
{\multiply\height by2
\puthmorphism(\xpos,\ypos)[#4`#5`{#7}]{\height}{\arrowtypea}a}%
\putmorphism(\xpos,\ypos)(1,-1)[`#6`{#8}]{\height}{\arrowtypeb}l%
\advance\xpos by\height
\advance\xpos by\height
\putmorphism(\xpos,\ypos)(-1,-1)[``{#9}]{\height}{\arrowtypec}r%
}}
 
\def\putVtriangle{\@ifnextchar <{\putVtrianglep}{\putVtrianglep
   <\arrowtypea`\arrowtypeb`\arrowtypec;\height>}}
\def\Vtriangle{\@ifnextchar <{\Vtrianglep}{\Vtrianglep
   <\arrowtypea`\arrowtypeb`\arrowtypec;\height>}}
                                               
\def\Vtrianglep<#1>[#2`#3`#4;#5`#6`#7]{{
\settriparms[#1]
\width=\height                                 
\diagram
\putVtrianglep<\arrowtypea`\arrowtypeb`
\arrowtypec;\height>
(0,0)[#2`#3`#4;#5`#6`{#7}]
\enddiagram
}}
 
\def\putVtrianglepairp<#1>(#2)[#3;#4`#5`#6`#7`#8]{{
\settripairparms[#1]%
\setpos(#2)%
\settokens`#3`%
\advance\ypos by\height
\putmorphism(\xpos,\ypos)(1,-1)[`\tokend`{#6}]{\height}{\arrowtypec}l%
\puthmorphism(\xpos,\ypos)[\tokena`\tokenb`{#4}]{\height}{\arrowtypea}a%
\advance\xpos by\height
\puthmorphism(\xpos,\ypos)[\phantom{\tokenb}`\tokenc`{#5}]%
{\height}{\arrowtypeb}a%
\putvmorphism(\xpos,\ypos)[``{#7}]{\height}{\arrowtyped}m%
\advance\xpos by\height
\putmorphism(\xpos,\ypos)(-1,-1)[``{#8}]{\height}{\arrowtypee}r%
}}
 
\def\putVtrianglepair{\@ifnextchar <{\putVtrianglepairp}{\putVtrianglepairp%
    <\arrowtypea`\arrowtypeb`\arrowtypec`\arrowtyped`\arrowtypee;\height>}}
\def\Vtrianglepair{\@ifnextchar <{\Vtrianglepairp}{\Vtrianglepairp%
    <\arrowtypea`\arrowtypeb`\arrowtypec`\arrowtyped`\arrowtypee;\height>}}
                                               
\def\Vtrianglepairp<#1>[#2;#3`#4`#5`#6`#7]{{
\settripairparms[#1]
\settokens`#2`
\diagram
\putVtrianglepairp                             
<\arrowtypea`\arrowtypeb`\arrowtypec`
\arrowtyped`\arrowtypee;\height>
(0,0)[{#2};#3`#4`#5`#6`{#7}]
\enddiagram
}}

\def\putCtrianglep<#1>(#2,#3)[#4`#5`#6;#7`#8`#9]{{%
\settriparms[#1]%
\xpos=#2 \ypos=#3
\advance\ypos by\height
\putmorphism(\xpos,\ypos)(1,-1)[``{#9}]{\height}{\arrowtypec}l%
\advance\xpos by\height
\advance\ypos by\height
\putmorphism(\xpos,\ypos)(-1,-1)[#4`#5`{#7}]{\height}{\arrowtypea}l%
{\multiply\height by 2
\putvmorphism(\xpos,\ypos)[`#6`{#8}]{\height}{\arrowtypeb}r}%
}}
 
\def\putCtriangle{\@ifnextchar <{\putCtrianglep}{\putCtrianglep
    <\arrowtypea`\arrowtypeb`\arrowtypec;\height>}}
\def\Ctriangle{\@ifnextchar <{\Ctrianglep}{\Ctrianglep
    <\arrowtypea`\arrowtypeb`\arrowtypec;\height>}}
                                             
\def\Ctrianglep<#1>[#2`#3`#4;#5`#6`#7]{{
\settriparms[#1]
\width=\height                               
\diagram
\putCtrianglep<\arrowtypea`\arrowtypeb`
\arrowtypec;\height>
(0,0)[#2`#3`#4;#5`#6`{#7}]
\enddiagram
}}                                           
                                             
\def\putDtrianglep<#1>(#2,#3)[#4`#5`#6;#7`#8`#9]{{%
\settriparms[#1]%
\xpos=#2 \ypos=#3
\advance\xpos by\height \advance\ypos by\height
\putmorphism(\xpos,\ypos)(-1,-1)[``{#9}]{\height}{\arrowtypec}r%
\advance\xpos by-\height \advance\ypos by\height
\putmorphism(\xpos,\ypos)(1,-1)[`#5`{#8}]{\height}{\arrowtypeb}r%
{\multiply\height by 2
\putvmorphism(\xpos,\ypos)[#4`#6`{#7}]{\height}{\arrowtypea}l}%
}}
 
\def\putDtriangle{\@ifnextchar <{\putDtrianglep}{\putDtrianglep
    <\arrowtypea`\arrowtypeb`\arrowtypec;\height>}}
\def\Dtriangle{\@ifnextchar <{\Dtrianglep}{\Dtrianglep
   <\arrowtypea`\arrowtypeb`\arrowtypec;\height>}}
                                            
\def\Dtrianglep<#1>[#2`#3`#4;#5`#6`#7]{{
\settriparms[#1]
\width=\height                              
\diagram
\putDtrianglep<\arrowtypea`\arrowtypeb`
\arrowtypec;\height>
(0,0)[#2`#3`#4;#5`#6`{#7}]
\enddiagram
}}                                          
                                            
\def\setrecparms[#1`#2]{\width=#1 \height=#2}%
 
\def\recursep<#1`#2>[#3;#4`#5`#6`#7`#8]{{\m@th
\width=#1 \height=#2
\settokens`#3`
\settowidth{\tempdimen}{$\tokena$}
\ifdim\tempdimen=0pt
  \savebox{\tempboxa}{\hbox{$\tokenb$}}%
  \savebox{\tempboxb}{\hbox{$\tokend$}}%
  \savebox{\tempboxc}{\hbox{$#6$}}%
\else
  \savebox{\tempboxa}{\hbox{$\hbox{$\tokena$}\times\hbox{$\tokenb$}$}}%
  \savebox{\tempboxb}{\hbox{$\hbox{$\tokena$}\times\hbox{$\tokend$}$}}%
  \savebox{\tempboxc}{\hbox{$\hbox{$\tokena$}\times\hbox{$#6$}$}}%
\fi
\ypos=\height
\divide\ypos by 2
\xpos=\ypos
\advance\xpos by \width
\bfig
\putCtrianglep<-1`1`1;\ypos>(0,0)[`\tokenc`;#5`#6`{#7}]%
\puthmorphism(\ypos,0)[\tokend`\usebox{\tempboxb}`{#8}]{\width}{-1}b%
\puthmorphism(\ypos,\height)[\tokenb`\usebox{\tempboxa}`{#4}]{\width}{-1}a%
\advance\ypos by \width
\putvmorphism(\ypos,\height)[``\usebox{\tempboxc}]{\height}1r%
\efig
}}
 
\def\recurse{\@ifnextchar <{\recursep}{\recursep<\width`\height>}}
 
\def\puttwohmorphisms(#1,#2)[#3`#4;#5`#6]#7#8#9{{%
\puthmorphism(#1,#2)[#3`#4`]{#7}0a
\ypos=#2
\advance\ypos by 20
\puthmorphism(#1,\ypos)[\phantom{#3}`\phantom{#4}`#5]{#7}{#8}a
\advance\ypos by -40
\puthmorphism(#1,\ypos)[\phantom{#3}`\phantom{#4}`#6]{#7}{#9}b
}}
 
\def\puttwovmorphisms(#1,#2)[#3`#4;#5`#6]#7#8#9{{%
\putvmorphism(#1,#2)[#3`#4`]{#7}0a
\xpos=#1
\advance\xpos by -20
\putvmorphism(\xpos,#2)[\phantom{#3}`\phantom{#4}`#5]{#7}{#8}l
\advance\xpos by 40
\putvmorphism(\xpos,#2)[\phantom{#3}`\phantom{#4}`#6]{#7}{#9}r
}}
 
\def\puthcoequalizer(#1)[#2`#3`#4;#5`#6`#7]#8#9{{%
\setpos(#1)%
\puttwohmorphisms(\xpos,\ypos)[#2`#3;#5`#6]{#8}11%
\advance\xpos by #8
\puthmorphism(\xpos,\ypos)[\phantom{#3}`#4`#7]{#8}1{#9}
}}
 
\def\putvcoequalizer(#1)[#2`#3`#4;#5`#6`#7]#8#9{{%

\setpos(#1)%
\puttwovmorphisms(\xpos,\ypos)[#2`#3;#5`#6]{#8}11%
\advance\ypos by -#8
\putvmorphism(\xpos,\ypos)[\phantom{#3}`#4`#7]{#8}1{#9}
}}
 
\def\putthreehmorphisms(#1)[#2`#3;#4`#5`#6]#7(#8)#9{{%

\setpos(#1) \settypes(#8)
\if a#9 %
     \vertsize{\tempcounta}{#5}%
     \vertsize{\tempcountb}{#6}%
     \ifnum \tempcounta<\tempcountb \tempcounta=\tempcountb \fi
\else
     \vertsize{\tempcounta}{#4}%
     \vertsize{\tempcountb}{#5}%
     \ifnum \tempcounta<\tempcountb \tempcounta=\tempcountb \fi
\fi
\advance \tempcounta by 60
\puthmorphism(\xpos,\ypos)[#2`#3`#5]{#7}{\arrowtypeb}{#9}
\advance\ypos by \tempcounta
\puthmorphism(\xpos,\ypos)[\phantom{#2}`\phantom{#3}`#4]{#7}{\arrowtypea}{#9}
\advance\ypos by -\tempcounta \advance\ypos by -\tempcounta
\puthmorphism(\xpos,\ypos)[\phantom{#2}`\phantom{#3}`#6]{#7}{\arrowtypec}{#9}
}}
 
\def\setarrowtoks[#1`#2`#3`#4`#5`#6]{%
\def\toka{#1}
\def\tokb{#2}
\def\tokc{#3}
\def\tokd{#4}
\def\toke{#5}
\def\tokf{#6}
}
\def\hex{\@ifnextchar <{\hexp}{\hexp<1000`400>}}
\def\hexp<#1`#2>[#3`#4`#5`#6`#7`#8;#9]{%
\setarrowtoks[#9]
\yext=#2 \advance \yext by #2
\xext=#1 \advance\xext by \yext
\bfig
\putCtriangle<-1`0`1;#2>(0,0)[`#5`;\tokb``\tokd]
\xext=#1 \yext=#2 \advance \yext by #2
\putsquare<1`0`0`1;\xext`\yext>(#2,0)[#3`#4`#7`#8;\toka```\tokf]
\advance \xext by #2
\putDtriangle<0`1`-1;#2>(\xext,0)[`#6`;`\tokc`\toke]
\efig
}

\makeatother

 \newcount\beziercnt
 
 \newcount\grcalca
 \newcount\grcalcb
 \newcount\grcalcc
 \newcount\grcalcd
 \newcount\mordiam
 \newcommand\thlines{\thinlines}

 \newcommand{\mult}{                       
   \thlines                                
   \put(5,5){\oval(10,10)[b]} }

 \newcommand{\bmult}[1]{                   
   \thlines                                
   \grcalca = #1
   \divide \grcalca by 2
   \put(\grcalca,5){\oval(#1,10)[b]} }

 \newcommand{\braid}{                       
   \thlines                                 
   \bezier{\beziercnt}(10,10)(10,7.5)(7,6)  
   \bezier{\beziercnt}(0,0)(0,2.5)(3,4)     
   \bezier{\beziercnt}(0,10)(0,7.5)(5,5)    
   \bezier{\beziercnt}(10,0)(10,2.5)(5,5) }

 \newcommand{\ibraid}{                      
   \thlines                                 
   \bezier{\beziercnt}(10,10)(10,7.5)(5,5)  
   \bezier{\beziercnt}(0,0)(0,2.5)(5,5)     
   \bezier{\beziercnt}(0,10)(0,7.5)(3,6)    
   \bezier{\beziercnt}(10,0)(10,2.5)(7,4) }

 \newcommand{\mor}[2]{                      
   \thlines                                  
   \grcalca = #2                            
   \advance \grcalca by -\mordiam          
   \divide \grcalca by 2                    
   \put(0,#2) {\line(0,-1){\grcalca}}       
   \put(0,0) {\line(0,1){\grcalca}}         
   \grcalca = #2                            
   \divide \grcalca by 2                    
   \put(0,\grcalca) {\circle{\mordiam}}     
   \put(0,\grcalca)                         
       {\makebox(0,0){$\scriptstyle #1$}} }

  \newcommand{\multmor}[2]{                    
    \thlines                                   
    \grcalca = #2                              
    \advance \grcalca by 5                     
    \put(-2.5,0){\framebox                     
             (\grcalca,10){$\scriptstyle #1$}}}

 \newcommand{\twist}[2]{                    
   \thlines                                 
   \grcalca = #1                            
   \divide \grcalca by 2                    
   \grcalcb = #2                            
   \divide \grcalcb by 2                    
   \grcalcc = #2                            
   \divide \grcalcc by 3                    
   \bezier{\beziercnt}(0,0)(0,\grcalcc)(\grcalca,\grcalcb) 
   \grcalcd = #2
   \advance \grcalcd by -\grcalcc
   \multiply \grcalcc by 3
   \bezier{\beziercnt}(\grcalca,\grcalcb)(#1,\grcalcd)(#1,#2) }

 \newcommand{\idgr}[1]{                     
   \thlines                                 
   \put(0,0) {\line(0,1){#1}}}

 \newcommand{\objo}[1]{
   \put(0,1){\makebox(0,0)[b]{$#1$}}}

 \newcommand{\obju}[1]{
   \put(0,-10){\makebox(0,0)[b]{$#1$}}}

 \newcommand{\bgr}[2]{\vskip2ex
   \begin{center}
   \unitlength=0.22ex
   \grcalca = #2
   \advance \grcalca by 7
   \begin{picture}(#1,\grcalca)}

 \newcommand{\egr}{
   \end{picture}{\vskip2.5ex}\end{center} }

\newtheorem{thm}{Theorem}[section]

\newtheorem{cor}[thm]{Corollary}
\newtheorem{lma}[thm]{Lemma}
\newtheorem{prop}[thm]{Proposition}

\newtheorem{defn}[thm]{Definition}
\newtheorem{eg}[thm]{Example}

\newcommand{\tensor}{\otimes}
\newcommand{\iso}{\cong}

\newcommand\Hom{{\mathop{\mathrm {Hom}}}} 
\newcommand\Ker{{\mathop{\mathrm {Ker}}}} 
 
\newcommand\rend{{\mathop{\mathrm {end}}}}
\newcommand\Aut{{\mathop{\mathrm {Aut}}}} 
\newcommand\Der{{\mathop{\mathrm {Der}}}}

\newcommand\Mod{{\mathop{\mathrm {\hbox{-}Mod}}}} 
\newcommand\comod{{\mathop{\mathrm {\hbox{-}comod}}}}

\newcommand\id{{\mathop{\mathrm {id}}}} 
\newcommand\ev{{\mathop{\mathrm {ev}}}}

\newcommand\og{{\mathop{\frak {g}}}} 
\newcommand\C{{\cal C}} 
\newcommand\M{{\cal M}}

\newcommand\YD{{\cal Y}{\cal D}}

\newcommand\x{{\mbox{-}}} 
\newcommand\X{{\mbox{--}}} 
 
\renewcommand\phi{\varphi}

\begin{document}

\title
 [On Lie Algebras in the Category of\\ Yetter-Drinfeld Modules] 
 {On Lie Algebras in the Category of\\ Yetter-Drinfeld Modules} 

\author{Bodo Pareigis}

\address{Mathematisches Institut der Universit\"at\\
 Theresienstr.39\\
 80333 M\"unchen\\ 
 Germany} 

\email{pareigis@@rz.mathematik.uni-muenchen.de}


\date{\today}

\maketitle


 \begin{abstract}
 The category of Yetter-Drinfeld modules $\YD^K_K$ over a 
Hopf algebra $K$ (with bijektive antipode over a field $k$) 
is a braided monoidal category. If $H$ is a Hopf algebra in 
this category then the primitive elements of $H$ do not 
form an ordinary Lie algebra anymore. We introduce the 
notion of a (generalized) Lie algebra in $\YD^K_K$ such 
that the set of primitive elements $P(H)$ is a Lie algebra 
in this sense. Also the Yetter-Drinfeld module of 
derivations of an algebra $A$ in $\YD^K_K$ is a Lie 
algebra. Furthermore for each Lie algebra in $\YD^K_K$ 
there is a universal enveloping algebra which turns out to 
be a Hopf algebra in $\YD^K_K$. 


 \end{abstract}

 \section{Introduction}

 The concept of Hopf algebras in braided categories has 
turned out to be very important in the context of 
understanding the structure of quantum groups and 
noncommutative noncocommutative Hopf algebras. In 
particular the work of Radford \cite{Ra}, Majid \cite{Ma}, 
Lusztig \cite{Lu}, and Sommerh\"auser \cite{So} show the 
importance of the decomposition of quantum groups into a 
product of ordinary Hopf algebras and of Hopf algebras in 
braided categories. 

 Since by the work of Yetter \cite{Ye} Hopf algebras in 
braided categories that are defined on an underlying 
(finite-dimensional) vector space can be considered as Hopf 
algebras in some category of Yetter-Drinfeld modules, we 
will restrict our attention to Hopf algebras $H$ in a 
category of Yetter-Drinfeld modules $\YD_K^K$ over a 
Hopf algebra $K$ with bijective antipode. 

 There are two structurally interesting and important 
concepts that survive in this generalized situation, the 
concept of group-like elements ($\Delta(g) = g \tensor g$, 
$\varepsilon(g) = 1$) and the concept of primitive elements 
($\Delta(x) = x \tensor 1 + 1 \tensor x$, $\varepsilon(x) = 
0$). 

 For ordinary Hopf algebras $H$ the set of primitive 
elements $P(H)$ of $H$ forms a Lie algebra. This result (in 
a somewhat generalized form) still holds for Hopf algebras 
in a symmetric monoidal category. This is, however, not 
true for braided monoidal categories.

 There have been various attempts to generalize the notion 
of Lie algebras to braided monoidal categories. The main 
obstruction for such a generalization is the assumption 
that the category is only braided and not symmetric. One of 
the most important examples of such braided categories is 
given by the category of Yetter-Drinfeld modules $\YD_K^K$ 
over a Hopf algebra $K$ with bijective antipode which is 
always properly braided (except for $K = k$, the base 
field) \cite{P2}. 

 We introduce a concept of Lie algebras in $\YD_K^K$ 
that generalizes the concepts of ordinary Lie algebras, Lie 
super algebras, Lie color algebras, and $(G,\chi)$-Lie 
algebras as given in \cite{P1}. 

 The Lie algebras defined on Yetter-Drinfeld modules have 
{\em partially defined $n$-ary} bracket operations for every $n 
\in {\Bbb N}$ and every primitive $n$-th root of unity. 
They satisfy generalizations of the (anti-)symmetry and 
Jacabi identities. 

 Our main aim is to show that these Lie algebras have 
universal enveloping algebras which turn out to be Hopf 
algebras in $\YD_K^K$. Conversely the set of primitive 
elements of a Hopf algebra in $\YD_K^K$ is such a 
generalized Lie algebra. We also give an example that 
generalizes the concept of orthogonal or symplectic Lie 
algebras. 

 \section{Braid Symmetrization} \label{braidsymm}

 We begin with two simple module theoretic observations. The 
following is well known: if $A, B$ are algebras and $M$ is 
an $A$-$B$-bimodule, then $\Hom_A(.P,.M)$ is a right $B$-module 
for every $A$-module $P$. We need a comodule analogue of 
this.

 Let $A$ be an algebra, $C$ be a coalgebra, and ${}_AM^C$ be 
an $A$-$C$-dimodule, i.e. a left $A$-module and a right
 $C$-comodule such that $\delta(am) = (a \tensor 
1)\delta(m)$. 

 \begin{prop} \label{induceddimodule}
 Let $P$ be a finitely generated left $A$-module. Then 
$\Hom_A(.P,.M)$ is a right $C$-comodule with the canonical 
comodule structure such that
 $$\big(\Hom_A(P,M) \buildrel \delta \over \longrightarrow 
\Hom_A(P,M) \tensor C \to \Hom_A(P,M \tensor C) \big) = 
\Hom_A(P,\delta).$$
 \end{prop}

 \begin{pf}
 Let $p_1, \ldots, p_n$ be a generating set of $P$ and let 
$f \in \Hom_A(.P,.M)$. Let $m_i := f(p_i)$. Then by the 
structure theorem on comodules the $m_i$ are contained in a 
finite dimensional subcomodule $M_0 \subseteq M$ which is 
even a comodule over a finite dimensional subcoalgebra $C_0 
\subseteq C$, i.e. the diagram
 $$\bfig
 \putmorphism(0, 400)(1, 0)[M_0`M_0 \tensor C_0`\delta]{500}1a
 \putmorphism(0, 0)(1, 0)[M`M \tensor C`\delta]{500}1a
 \putmorphism(0, 400)(0, -1)[``]{400}1l
 \putmorphism(500, 400)(0, -1)[``]{400}1r
 \efig$$
 commutes. Furthermore $M_1 := AM_0$ is a $C_0$-comodule 
contained in $M$, since $M$ is a dimodule, and $f: P \to M$ 
obviously factors through $M_1$. Since $M$ and $M_1$ are 
dimodules the diagram 
 $$\bfig
 \putmorphism(0, 400)(1, 0)[\Hom_A(.P,.M_1)`\Hom_A(.P,.M_1 
\tensor C_0) `\delta_*]{1100}1a
 \putmorphism(0,   0)(1, 0)[\Hom_A(.P,.M)  `\Hom_A(.P,.M 
\tensor C  ) `\delta_*]{1100}1a
 \putmorphism(1100, 400)(1, 0)[\phantom{\Hom_A(.P,.M_1 
\tensor C_0)}`\Hom_A(.P,.M_1) \tensor C_0`\iso]{1100}{-1}a 
 \putmorphism(1100,   0)(1, 0)[\phantom{\Hom_A(.P,.M
\tensor C  )}`\Hom_A(.P,.M) \tensor C`]{1100}{-1}a
 \putmorphism(0, 400)(0, -1)[``]{400}1r
 \putmorphism(1100, 400)(0, -1)[``]{400}1r
 \putmorphism(2200, 400)(0, -1)[``]{400}1r
 \efig$$
 commutes, so each $f$ has a uniquely defined image 
$\delta_*(f) \in \Hom_A(.P,.M) \tensor C  $. Now it is easy 
to check that this map induces a comodule structure on 
$\Hom_A(P,M)$. 
 \end{pf}

 The second observation is the following. We consider
 $k$-algebras $A$ and $B$. Let $\alpha: B \to A$ be 
an algebra homomorphism. $\alpha$ induces an underlying 
functor $V_\alpha: A\Mod \to B\Mod$ with right adjoint 
$\Hom_B(A,\x): B\Mod \to A\Mod$. If $\alpha: B \to A$ is 
surjective then $\Hom_B(A,M) \to \Hom_B(B,M) \iso M$ is 
injective, so that we can identify $\Hom_B(A,M) = \{m \in M 
| \Ker(\alpha)m = 0\}$.  

 Let $B_n$ be the Artin braid group with generators 
$\tau_i$, $i = 1, \ldots,n-1$ and relations 
 \begin{equation}\begin{array}{rll}
 \tau_i\tau_j &\hskip-1.2ex\  = \tau_j\tau_i & \mbox{ if } 
|i-j| \geq 2; \\ 
 \tau_i\tau_{i+1}\tau_i &\hskip-1.2ex\ = 
\tau_{i+1}\tau_i\tau_{i+1}.&\\ 
 \end{array}\end{equation}
 Let $\zeta \in k$ be invertible. Then $kB_n \ni \tau_i 
\mapsto \zeta \tau_i \in kB_n$ (for the generators $\tau_i$ 
of $B_n$) is an algebra automorphism denoted again by 
$\zeta: kB_n \to kB_n$. This holds true since the relations 
for $B_n$ are homogeneous. 

 (Observe that this construction can be performed for every 
group algebra if the group is given by generators and 
homogeneous relations. The given construction of an 
automorphism for every  $\zeta \in U(k)$ defines a group 
homomorphism $U(k) \to \Aut(kB_n) \to \Aut(kB_n\Mod)$.) 

 Now consider the canonical quotient homomorphism $B_n \to 
S_n$ from the braid group onto the symmetric group. It 
induces a surjective homomorphism $\gamma: kB_n \to kS_n$ 
with kernel 
 $$\Ker(\gamma) = \langle \psi(\tau_i^2 - 1)\phi | \phi, 
\psi \in B_n, i = 1, \ldots,n-1 \rangle.$$

 The composition $\alpha: kB_n \buildrel \zeta \over 
\longrightarrow kB_n \buildrel \gamma \over \longrightarrow 
kS_n$ defines a functor $kB_n\Mod \to kS_n\Mod$ by

 \begin{equation}\begin{array}{rl}
 M(\zeta) :&= \Hom_{kB_n}({}_{\zeta}kS_n,M) \\
 &= \{m \in M|\phi^{-1}\tau_i^2\phi(m) = \zeta^2m\quad 
\forall \phi \in B_n,\ i = 1, \ldots, n-1\}. \\
 &= \{m \in M|\tau_i^2\phi(m) = \zeta^2\phi(m)\quad 
\forall \phi \in B_n,\ i = 1, \ldots, n-1\}. \\
 \end{array}
 \end{equation}
 This holds since the map $\gamma: kB_n \to kS_n$ has as 
kernel the two-sided ideal generated as a $k$-subspace by 
$\{\psi (\tau_i^2 - 1) \phi\ | \psi, \phi \in B_n, i = 1, 
\ldots, n-1 \}$. So $f \in \Hom_{kB_n}({}_{ \zeta}kS_n,M)$ 
with $f(1) = m \in M$, iff $\zeta^{-1}(\psi (\tau_i^2 - 1) 
\phi)m = 0$ for all $\psi, \phi, i$, iff $\zeta^{-1}(\tau_i^2 - 
1) \phi m = 0$ for all $\phi, i$, iff $\tau_i^2 \phi m = 
\zeta^2\phi m$ for all $\phi, i$, iff 
 $\phi^{-1}\tau_i^2\phi(m) = \zeta^2m$ for all $\phi, i$. 

 If the action of $B_n$ on $M$ is given by an action of 
$S_n$ and the canonical epimorphism $B_n \to S_n$, then the 
construction of the $M(\zeta)$ becomes trivial, since 
$M(\zeta) = \{m \in M | \tau_i^2\phi(m) = \zeta^2 \phi(m) = 
\phi(m) \} = 0$ if $\zeta^2 \not= 1$ and $M(-1) = M(1) = 
M$. Observe that the module $M(\zeta)$ depends only on 
$\zeta^2$, but that the action of $kS_n$ on $M(\zeta)$ 
depends on $\zeta$.

 $M(1)$ gives a solution of the following universal 
problem. 

 \begin{prop}
 For every $kB_n$-module $M$ the subspace
 $$M(1) := \{m \in M|\phi^{-1}\tau_i^2\phi(m) = m\quad 
\forall \phi \in B_n,\ i = 1, \ldots, n-1\}$$
 is a $kS_n$-module and the inclusion $M(1) \to M$ is a 
$kB_n$-module homomorphism, such that for every
 $kS_n$-module $T$ and every $kB_n$-module homomorphism $f:T 
\to M$ there is a unique $kS_n$-module homomorphism $g: T 
\to M(1)$ such that the diagram 
  $$\bfig
 \putmorphism(0, 0)(1, 0)[M(1)`M`\iota]{800}1a
 \putmorphism(0, 400)(0, 1)[T``g]{400}1l
 \putmorphism(0, 400)(2, -1)[``f]{800}1r
 \efig$$
commutes.
 \end{prop}

 \begin{defn}{\em 
 We call the functor $\x(\zeta): kB_n\Mod \to kS_n\Mod$ the 
$\zeta$-{\em sym\-me\-tri\-za\-tion} of $kB_n$-modules.} 
 \end{defn}

 The definition gives 
 $$M(\zeta) = \{m \in M| \phi^{-1} \tau_i^2 \phi(m) = 
\zeta^2 m\quad \forall \phi \in B_n,\ i = 1, \ldots, n-
1\}.$$ 
 The action of $S_n$ on $M(\zeta)$ is given by 
 \begin{equation}\label{zetasymm}
 \sigma_i(m) = \zeta^{-1}\tau_i(m),
 \end{equation}
 where $\sigma_i$ resp. $\tau_i$ are the canonical 
generators of $S_n$ resp. $B_n$. Thus $M(\zeta)$ is also a 
$kB_n$-submodule of $M$. Since the functor $M \mapsto 
M(\zeta)$ is a rightadjoint functor, it preserves limits. 
Like for eigenspaces we have that the sum of the subspaces 
$M(\zeta)$ for all $\zeta$ with different $\zeta^2$ is a 
direct sum. On $M(\zeta)$ we have two distinct
 $kS_n$-structures $\sigma_i(m) = -\zeta^{-1}\tau_i(m)$ and 
$\sigma_i = \zeta^{-1}\tau_i(m)$, since $\zeta$ and $-
\zeta$ define the same subspace $M(\zeta) = M(-\zeta) 
\subseteq M$. 

The $\zeta$-symmetrization $M(\zeta)$ of $M$ can also be 
calculated by 

 \begin{lma} \label{append1}
 $$M(\zeta) = \{m \in M | \tau_i^{- 1} \tau_{i+1}^{- 1} 
\ldots \tau_{j-1}^{- 1} \tau_{j}^2 \tau_{j-1} \ldots 
\tau_{i+1} \tau_i(m) = \zeta^2 m\quad  \forall 1 
\leq i \leq j \leq n-1 \}$$ 
 which reduces the number of conditions to be imposed on the 
$m \in M$ in order to be in $M(\zeta)$.
 \end{lma} 

 \begin{pf}
 Given in  Appendix.
 \end{pf}

 One of the interesting $kB_n$-structures, for which we will 
apply the previous construction, occurs on $n$-fold tensor 
products $M^n := M \tensor \ldots \tensor M$ of an object 
$M$ in a braided monoidal category of vector spaces. 

 Let $K$ be a Hopf algebra. Let $M$ be an $K$-module such 
that $M$ is a $kB_n$-$K$-bimodule. The functoriality of our 
construction then makes $M(\zeta)$ again an $K$-module and 
in fact a $kS_n$-$K$-bimodule. 

 Let $M$ be an $K$-comodule such that $M$ is a
 $kB_n$-$K$-dimodule. Then by Proposition 
\ref{induceddimodule} $M(\zeta)$ is an $K$-comodule and in 
fact a $kS_n$-$K$-dimodule. 

 Let $K$ be a Hopf algebra with bijective antipode. Let 
$\YD_K^K$ denote the category of Yetter-Drinfeld modules 
over $K$, i.e. of right $K$-modules and right $K$-comodules 
$M$ such that $\sum (x \cdot h)_0 \tensor (x \cdot h)_1 = 
\sum (x_0 \cdot h_2) \tensor S(h_1)x_1h_3$ for all $x \in 
M$. The usual tensor product makes $\YD_K^K$ a monoidal 
category. $\YD_K^K$ has a braiding given by $\tau_{X,Y}: X 
\tensor Y \to Y \tensor X$, $\tau(x \tensor y) = \sum 
y_{0} \tensor xy_{1}$. We assume that the reader is 
familiar with the properties of the $B_n$-action that is 
induced by the braiding $\tau$ on $n$-fold tensor products 
(\cite{M} 10.6). 

 \begin{thm}\label{symmetrization}
 Let $K$ be a Hopf algebra with bijective antipode. Then 
for each $\zeta \in k^*$ and each $n \geq 2$ the 
construction given above defines a (non-additive) functor 
 $$\YD_K^K \ni M \mapsto (M \tensor \ldots \tensor 
M)(\zeta) \in \YD_K^K.$$
 \end{thm}

 \begin{pf}
 If $M \in \YD_K^K$ then the $n$-fold tensor product $M 
\tensor  \ldots \tensor M$ is a Yetter-Drinfeld module on 
which $B_n$ and thus $kB_n$ acts in such a way, that $M$ is 
a $(kB_n,K)$-bimodule and a $(kB_n,K)$-dimodule. The 
$\zeta$-symmetrization functor $\x(\zeta)$ preserves the 
module and comodule structures hence the Yetter-Drinfeld 
structure. 

 The functor is not additive since the ``diagonal'' functor 
$M \mapsto M \tensor M$ is not additive. 
 \end{pf}

 We abbreviate $M^n(\zeta) := M^{\tensor n}(\zeta) = M 
\tensor \ldots \tensor M(\zeta)$. Then $M^n(\zeta)$ is a 
submodule of $M^n$ in the category of Yetter-Drinfeld 
modules and the elements in $M^n(\zeta)$ are of the form $z 
= \sum_k x_{k,1} \tensor \ldots \tensor x_{k,n}$. We often 
suppress the summation index and summation sign and simply 
write $z = z_1 \tensor \ldots \tensor z_n \in M^n(\zeta)$ 
although $M^n(\zeta)$ does not decompose into a tensor 
product.

 \section{Symmetric Multiplication and Jacobi Identities}

 For the rest of the paper let $\C$ be the category 
$\YD_K^K$ of Yetter-Drinfeld modules over a Hopf algebra 
$K$ with bijective antipode over a field $k$. We study 
objects $P \in \C$ together with (partially defined) 
operations in $\C$ 
  $$[.,.]: P\tensor \ldots \tensor P(\zeta) = P^n(\zeta) 
\to P$$ 
 for all $n \in \Bbb N$ and all {\em primitive $n$-th roots 
of unity} $\zeta$. 

 Occasionally we write $[.,.]_n$ for such an operation 
$[.,.]$. By composing such operations certain additional 
operations may be constructed as follows. 

 \begin{prop} \label{jacobidef1}
 Let $\zeta$ be a primitive $n$-th roots of unity. Then the 
operations 
 $$[.,[.,.]_n]_2: P^{n+1}(\zeta) \ni x_1 \tensor \ldots \tensor 
x_{n+1} \mapsto [x_1,[x_2, \ldots,x_{n+1}]] \in P$$ 
 and 
 $$[[.,.]_n,.]_2: P^{n+1}(\zeta) \ni x_1 \tensor \ldots \tensor 
x_{n+1} \mapsto [[x_1, \ldots,x_n],x_{n+1}] \in P$$ 
 are well defined. 
 \end{prop} 

 \begin{pf}
 Given in  Appendix.
 \end{pf}

 We will have to consider objects
 $$P^{n+1}(-1,\zeta) := P \tensor P^n(\zeta) \cap \{z \in 
P^{n+1} | \forall \phi \in S_n : (1 \tensor \phi)^{-1} 
\tau_1^2 (1 \tensor \phi)(z) = z \}.$$ 
 Since this is a kernel (limit) construction in $\C$, 
$P^{n+1}(-1,\zeta)$ is again an object in $\C$. 

 \begin{prop} \label{jacobidef2}
 Let $\zeta$ be a primitive $n$-th roots of unity. Then the 
operations 
 $$[.,[.,.]_n]_2: P^{n+1}(-1,\zeta) \ni x \tensor y_1 
\tensor \ldots \tensor y_n \mapsto [x,[y_1, \ldots, y_n]] 
\in P$$                            
 and 
 $$[.,[.,.]_2,.]_n \tau_{i-1} \ldots \tau_1: P^{n+1}(-
1,\zeta) \ni x \tensor y_1 \tensor \ldots \tensor y_n 
\mapsto [y_1, \ldots , [x, y_i], \ldots, y_n] \in P$$ 
 are well defined. 
 \end{prop} 

 \begin{pf}
 Given in  Appendix.
 \end{pf}

 We introduce special bracket multiplications which then 
lead to the definition of a Lie algebra on a
 Yetter-Drinfeld module. 

 \begin{defn} \label{symm_mult}{\em
 Let $A$ be an algebra in $\C = \YD_K^K$ and let $\nabla^n: A \tensor 
\ldots \tensor A \to A$ denote the $n$-fold multiplication. 
We define a {\em bracket} or {\em symmetric multiplication} 
 $$[.,.]: A^n(\zeta) \to A \qquad \mbox{ by } \qquad [z] := 
\sum_{\sigma \in S_n} \nabla^n \sigma(z)$$ 
 where the action of $S_n$ on $A^n(\zeta)$ is given as in 
(\ref{zetasymm}).}
 \end{defn}

 We will only use those bracket operations which are 
defined with $\zeta$ a primitive $n$-th root of unity (for 
all $n \in \Bbb N$ and all $\zeta$).

 We consider these bracket operations as a generalization of 
the Lie-bracket $[\x,\x]: L \times L \to L$ or $[\X]: L 
\tensor L \to L$. Observe that our bracket operation is only 
partially defined and should {\em not} be considered as a 
multilinear operation, since $A^n(\zeta) \subseteq A^n$ is 
just a submodule in $\C$ and does not necessarily decompose 
into an $n$-fold tensor product. The elements in 
$A^n(\zeta)$ are, however, of the form $z = \sum_k x_{k,1} 
\tensor \ldots \tensor x_{k,n}$. 

 If we suppress the summation index and the summation sign 
then we may write the bracket operation on $z = x_1 \tensor 
\ldots \tensor x_n$ also as $[z] = [x_1, \ldots, x_n]$. If 
we define 
 \begin{equation}\label{sigma_action}
 \sigma(z) =: x_{\sigma^{-1}(1)} \tensor \ldots \tensor 
x_{\sigma^{-1}(n)}
 \end{equation}
 then we get
 $$[x_1, \ldots, x_n] = \sum_{\sigma \in S_n} 
x_{\sigma(1)}\cdot \ldots \cdot x_{\sigma(n)}.$$
 Observe that the components $x_1, \ldots, x_n$ in this 
expression are interchanged according to the action of the 
braid group resp. the symmetric group on $A^n(\zeta)$, so 
$x_{\sigma(1)} \tensor \ldots \tensor x_{\sigma(n)}$ is 
only a symbolic expression. 

 The bracket operation obviously satisfies the
 ``anti''-symmetry identity 
 \begin{equation}\label{symm}
 [\sigma(z)] = [z] \qquad \qquad \forall \sigma \in S_n.
 \end{equation}

 We apply Proposition \ref{jacobidef1} to an algebra $A$ in 
$\C$ with the operations given in Definition 
\ref{symm_mult} and get 

 \begin{thm} {\em (1.~Jacobi identity)} \label{1.jac}
 For all $n \in {\Bbb N}$, for all primitive $n$-th roots of 
unity $\zeta$, and for all $z \in A^{n+1}(\zeta)$ we have
 \begin{equation}\label{Jacobi1}
 \sum_{i=1}^{n+1} [x_i,[x_1,\ldots,\hat x_i, \ldots, x_{n+1}]] = 
 \sum_{i=1}^{n+1} [.,[.,.]_n]_2(1 \ldots i)(z) = 0,
 \end{equation} 
 where we use the notation {\em (\ref{sigma_action})}.
 \end{thm}

 \begin{pf}
 We define $(S_{n+1})_{(i)} := \{ \sigma \in S_{n+1} | 
\sigma(i) = 1\}$. Then $S_{n+1} = \bigcup_i 
(S_{n+1})_{(i)} $. For $\sigma \in (S_{n+1})_{(i)}$ let 
$\bar{\rho} := \sigma(i \ldots 1)$. Since $\bar{\rho}(1) = 
1$ there is a unique $\rho \in S_n$ with $\bar{\rho} = 1 
\tensor \rho$ and $\sigma = (1 \tensor \rho)(1 \ldots i)$. 
So we obtain a bijection
 $$S_n \ni \rho \mapsto (1 \tensor \rho)(1 \ldots i) \in 
(S_{n+1})_{(i)}.$$
 Analogously we define $(S_{n+1})^{(i)} := \{ \sigma \in S_{n+1} | 
\sigma(i) = n+1\}$ and get a bijection
 $$S_n \ni \rho \mapsto (\rho \tensor 1)(n+1 \ldots i) \in 
(S_{n+1})^{(i)}.$$
 Now observe that $\tau_n \ldots \tau_1(z) = \zeta^n 
\sigma_n \ldots \sigma_1(z) = (n+1 \ldots 1)(z)$ (by 
$\zeta^n = 1$) for $z \in P^{n+1}(\zeta)$ to get 
 $$\begin{array}{l @{\ } l}
 \sum_{i=1}^n &[.,[.,.]](1 \ldots i)(z) =\\
 &= \sum_{i=1}^n \nabla(1 \tensor [.,.]) (1 \ldots i)(z) - 
\nabla([.,.] \tensor 1)\tau_n \ldots \tau_1(1 \ldots i)(z)\\ 
 &= \sum_{i=1}^n \nabla(1 \tensor [.,.]) (1 \ldots i)(z) - 
\nabla([.,.] \tensor 1)(n+1 \ldots i)(z)\\ 
 &= \sum_{i=1}^n \sum_{\rho \in S_n} \nabla^{n+1}(1 \tensor 
\rho)(1 \ldots i)(z) - \nabla^{n+1} (\rho \tensor 1)(n+1 
\ldots i)(z)\\ 
 &= \sum_{\sigma \in S_{n+1}} \nabla^{n+1} \sigma(z) - 
\nabla^{n+1} \sigma(z) = 0.\\ 
 \end{array}$$
 \end{pf}

 \begin{thm} {\em (2.~Jacobi identity)} \label{2.jac}
 For all $n \in {\Bbb N}$, for all primitive $n$-th roots 
of unity $\zeta$, and for all $z = x \tensor y_1 \tensor 
\ldots \tensor y_n \in A^{n+1}(-1,\zeta)$ we have 
  \begin{equation}\label{Jacobi2}
 [x,[y_1, \ldots, y_n]] = \sum_{i=1}^n [y_1, \ldots , [x, 
y_i], \ldots, y_n]
  \end{equation}
 where $y_1 \tensor \ldots \tensor y_{i-1} \tensor x 
\tensor y_i \tensor \ldots \tensor y_n := \tau_{i-1} \ldots 
\tau_1(z)$ and 
 \begin{equation} \label{sigma_action2}
 [y_1, \ldots, y_{i-1}, 
[x, y_i], \ldots, y_n] = [.,[.,.]_2,.]_n 
\tau_{i-1} \ldots \tau_1(z).
 \end{equation} 
 \end{thm}

 \begin{pf}
 The equation in the Theorem can also be written as
 $$[.,[.,.]_n]_2(z) = \sum_{i=1}^n [.,[.,.]_2,.]_n \tau_{i-
1} \ldots \tau_1(z).$$ 

 Lemma \ref{equn5} together with $\tilde{\phi}(i) = j$ 
shows 
 $$\begin{array}{r @{\ } l}
 \nabla^n\phi(1 \tensor \ldots \tensor \nabla &\tensor 
\ldots \tensor 1)\tau_{i-1} \ldots \tau_1(z)\\
 &= \nabla^n (1 \tensor \ldots \tensor \nabla \tensor 
\ldots \tensor 1) \phi_{(i)} \tau_{i-1} \ldots \tau_1(z)\\
 &= \nabla^{n+1} \tau_{j-1} \ldots \tau_1 (1 \tensor 
\phi)(z);\\
 \nabla^n\phi(1 \tensor \ldots \tensor \nabla &\tensor 
\ldots \tensor 1)\tau_i \ldots \tau_1(z)\\
 &= \nabla^n (1 \tensor \ldots \tensor \nabla \tensor 
\ldots \tensor 1) \phi_{(i)} \tau_i \ldots \tau_1(z)\\
 &= \nabla^{n+1} \tau_j \ldots \tau_1 (1 \tensor 
\phi)(z);\\
 \end{array}$$
 hence 
 \begin{equation} \label{equn6}
 \begin{array}{r @{\ } l}
 \nabla^n\phi(1 \tensor \ldots \tensor \nabla &\tensor 
\ldots \tensor 1)\tau_{i-1} \ldots \tau_1(z)\\
 &= \nabla^{n+1} \tau_{k-1} \ldots \tau_1 (1 \tensor 
\phi)(z)\\
 &= \nabla^{n+1} \tau_l \ldots \tau_1 (1 \tensor 
\phi)(z)\\
 &= \nabla^n\phi(1 \tensor \ldots \tensor \nabla \tensor 
\ldots \tensor 1)\tau_j \ldots \tau_1(z).\\
 \end{array}
 \end{equation}
 for all $i,j = 1, \ldots, n$ with $\tilde{\phi}(i) = 
\tilde{\phi}(j) + 1$, i.e. for all $i$ except 
$\tilde{\phi}^{-1}(1)$ and all $j$ except $\tilde{\phi}^{-
1}(n)$. The other $i$'s and $j$'s used in (\ref{equn6}) are 
in bijective correspondence. 
 
 To prove the equation of the theorem we write each $\sigma 
\in S_n$ as $\zeta^r \phi$ with a representative $\phi \in 
B_n$ and a suitable power $\zeta^r$ according to 
(\ref{zetasymm}) and use (\ref{equn6}). Then we get
 $$\begin{array}{r @{} l}
 \sum_{i=1}^n &[.,[.,.]_2,.]_n \tau_{i-
1} \ldots \tau_1(z) = \\
 &= \sum_{i=1}^n [.,.]_n (1 \tensor \ldots \tensor (\nabla - 
\nabla \tau) \tensor \ldots \tensor 1) \tau_{i-1} \ldots 
\tau_1(z) = \\ 
 &= \sum_{i=1}^n \sum_{\sigma \in S_n} \nabla^n \zeta^r 
\phi (1 \tensor \ldots \tensor \nabla \tensor \ldots 
\tensor 1) \tau_{i-1} \ldots \tau_1(z) \\ 
 &\ \ -\sum_{j=1}^n \sum_{\sigma \in S_n} \nabla^n \zeta^r 
\phi (1 \tensor \ldots \tensor \nabla \tensor \ldots 
\tensor 1) \tau_j \tau_{j-1} \ldots \tau_1(z) \\ 
 &= \sum_{\sigma \in S_n} \nabla^{n+1} (1 \tensor \zeta^r \phi)(z) 
 - \sum_{\sigma \in S_n} \nabla^{n+1} \tau_n \ldots \tau_1 
(1 \tensor \zeta^r \phi)(z) \\
 &= \sum_{\sigma \in S_n} \nabla^{n+1} (1 \tensor \zeta^r \phi)(z) 
 - \sum_{\sigma \in S_n} \nabla^{n+1} (\zeta^r \phi \tensor 
1)\tau_n \ldots \tau_1 (z) \\ 
 &= \nabla (1 \tensor \sum_{\sigma \in S_n} \nabla^n 
\sigma)(z)
 - \nabla (\sum_{\sigma \in S_n} \nabla^n 
\sigma \tensor 1 )\tau_{P,P^n}(z)\\
 &= (\nabla - \nabla \tau)(1 \tensor [.,.]_n)(z)\\
 &= [.,[.,.]_n]_2(z).
 \end{array}$$
 \end{pf}

 Clearly there are symmetric right sided identities.

 \section{Lie Algebras on Yetter-Drinfeld Modules}

 Now we can define the notion of a Lie algebra in the category of 
Yetter-Drinfeld modules.

 \begin{defn}{\em 
 A Yetter-Drinfeld module $P$ together with operations in $\YD_K^K$
  $$[.,.]: P\tensor \ldots \tensor P(\zeta) = P^n(\zeta) \to P$$
 for all $n \in \Bbb N$ and all primitive $n$-th roots of unity 
$\zeta$ is called a {\em Lie algebra} if the following 
identities hold: 
 \begin{enumerate} 
 \item for all $n \in \Bbb N$, for all primitive $n$-th 
roots of unity $\zeta$, for all $\sigma \in S_n$, and for 
all $z \in P^n(\zeta)$ 
 $$[z] = [\sigma(z)],$$ 
 \item 
 for all $n \in {\Bbb N}$, for all primitive $n$-th roots of 
unity $\zeta$, and for all $z \in P^{n+1}(\zeta)$ 
 $$\sum_{i=1}^{n+1} [x_i,[x_1,\ldots,\hat x_i, \ldots, x_{n+1}]] = 
 \sum_{i=1}^{n+1} [.,[.,.]](1 \ldots i)(z) = 0,$$ 
 where we use the notation (\ref{sigma_action}),
 \item 
 for all $n \in {\Bbb N}$, for all primitive $n$-th roots 
of unity $\zeta$, and for all $z = x \tensor y_1 \tensor 
\ldots \tensor y_n \in P^{n+1}(-1,\zeta)$ we have 
 $$[x,[y_1, \ldots, y_n]] = \sum_{i=1}^n [y_1, \ldots , [x, 
y_i], \ldots, y_n]$$ 
 where we use the notation (\ref{sigma_action2}).
 \end{enumerate}} 
 \end{defn}

 \begin{cor} \label{algebras_are_Lie}
 Let $A$ be an algebra in $\YD_K^K$. Then $A$ 
carries the structure of a Lie algebra $A^L$ with 
the symmetric multiplications
 $$[\X]: A^n(\zeta) \to A \qquad \mbox{ by } \qquad [z] := 
\sum_{\sigma \in S_n} \nabla^n \sigma(z).$$ 
 for all $n \in \Bbb N$ and all roots of unity $\zeta \in k^*$. 
 \end{cor} 

 \begin{pf} 
 This is a rephrasing of the ``anti''-symmetry identity 
(\ref{symm}) and the Jacobi identities (\ref{Jacobi1}) and 
(\ref{Jacobi2}) in Theorems  \ref{1.jac} and \ref{2.jac}. 
 \end{pf} 

 \section{The Lie Algebra of Primitive Elements}

 Let $A$ be an algebra in $\C = \YD_K^K$. Then $A \tensor 
A$ is an algebra with the multiplication $A \tensor A 
\tensor A \tensor A \buildrel 1 \tensor \tau \tensor 1 
\over \longrightarrow A \tensor A \tensor A \tensor A 
\buildrel \nabla \tensor \nabla \over \longrightarrow A 
\tensor A$. Let $p: A \to A \tensor A$ be the map $p(x) := 
x \tensor 1 + 1 \tensor x$. Then $p (= 1 \tensor \eta + 
\eta \tensor 1)$ is in $\C$ but $p$ is not an algebra 
morphism. Let $p^n: A^n \to (A \tensor A)^n$ be the
 $n$-fold tensor product of $p$ with itself. 

 \begin{lma} \label{primYD}
 Let $H$ be a Hopf algebra in $\C$. Then $P(H) := \{x \in H 
| \Delta(x) = x \tensor 1 + 1 \tensor x\}$ is a 
Yetter-Drinfeld submodule of $H$ in $\C$. 
 \end{lma}

 \begin{pf}
 $P(H) = \Ker(\Delta - p)$.
 \end{pf}

 In particular we have $\delta(x) \in P(H) \tensor K$ and 
$x \lambda \in P(H)$ for all $x \in P(H)$ and all $\lambda 
\in K$. 

 \begin{lma}
 $p^n(A^n(\zeta)) \subseteq (A \tensor A)^n(\zeta)$.
 \end{lma}

 \begin{pf}
 By Theorem \ref{symmetrization} $p: A \to A \tensor A$ 
induces $p^n: A^n(\zeta) \to (A \tensor A)^n(\zeta)$.
 \end{pf}

 \begin{thm}\label{mainthm}
 Let $\zeta$ be a primitive $n$-th root of unity and let $z 
\in A^n(\zeta)$. Then
 $$[p^n(z)] = p([z]).$$
 \end{thm}

 \begin{pf}
 If $z = \sum_k x_{k,1} \tensor \ldots \tensor x_{k,n} \in 
A^n(\zeta)$ then the equation of the theorem reads as
 \begin{equation}
 \begin{array}{c}
  \left[ \sum_k (x_{k,1} \tensor 1 + 1 \tensor x_{k,1}) 
\tensor \ldots \tensor (x_{k,n} \tensor 1 + 1 \tensor 
x_{k,n}) \right]  = \\ 
 \left[ \sum_k x_{k,1} \tensor \ldots \tensor x_{k,n} 
\right] \tensor 1 + 1 \tensor \left[ \sum_k x_{k,1} \tensor \ldots 
\tensor x_{k,n} \right] .\\ 
 \end{array}
 \end{equation}
 We want to evaluate 
 $$\begin{array}{l} \left[ \sum_k (x_{k,1} \tensor 1 + 1 
\tensor x_{k,1}) \tensor \ldots \tensor (x_{k,n} \tensor 1 + 
1 \tensor x_{k,n}) \right]  \\ 
 \qquad = \sum_{\sigma \in S_n} \nabla^n \sigma (\sum_i (x_{k,1} 
\tensor 1 + 1 \tensor x_{k,1}) \tensor \ldots \tensor 
(x_{k,n} \tensor 1 + 1 \tensor x_{k,n}))\\
 \end{array} $$
 where $\sigma \in S_n$ operates on $p^n(z) \in (A \tensor 
A)^n(\zeta)$ as described in section \ref{braidsymm}.

 Let $x,y \in A$ and $\tau(x \tensor y) = \sum_i u_i \tensor 
v_i$. Then $(1 \tensor x)\cdot(y \tensor 1) = (\nabla 
\tensor \nabla) (\sum_i 1 \tensor u_i \tensor v_i \tensor 1) 
= \sum_i u_i \tensor v_i = \tau(x \tensor y) = \sum_i (u_i 
\tensor 1)\cdot(1 \tensor v_i)$. So we have 
 \begin{equation}\label{multipl}
 \begin{array}{l} 
 (x \tensor 1)(y \tensor 1) = (xy \tensor 1),\\ 
 (x \tensor 1)(1 \tensor y) = (x \tensor y),\\ 
 (1 \tensor x)(1 \tensor y) = (1 \tensor xy),\\ 
 (1 \tensor x)(y \tensor 1) = \tau(x \tensor y). 
 \end{array}
 \end{equation} 

 We expand a product $(x_1 \tensor 1 + 1 \tensor x_1) 
\cdot \ldots \cdot (x_n \tensor 1 + 1 \tensor x_n)$. It 
produces after multiplication $2^n$ summands, each a product 
of $n$ terms. A typical product is $(x_1 \tensor 1)(1 
\tensor x_2)(x_3 \tensor 1)\ldots $, some of the factors 
being of the form $x_j \tensor 1$, the others of the form $1 
\tensor x_j$. To evaluate such a product we use the rule of 
multiplication in $A \tensor A$ given by $\nabla_{A \tensor 
A} = (\nabla \tensor \nabla)(1 \tensor \tau \tensor 1)$. 

 To explain the following calculation we consider as an 
example the product $(x_1 \tensor 1)(1 \tensor x_2)(x_3 
\tensor 1)(1 \tensor x_4)(x_5 \tensor 1)$. It is calculated 
with the following braid diagram 
 \bgr{40}{30} 
 \put(10,20){\braid}
 \put(20,10){\braid}
 \put(30,20){\braid}
 \put(0,10){\idgr{20}} 
 \put(40,10){\idgr{10}} 
 \put(10,0){\idgr{20}} 
 \put(35,0){\idgr{5}} 
 \put(0,5){\bmult{20}} 
 \put(30,5){\mult} 
 \egr 
 The second and fourth factors are pulled over to the right 
and then all factors are multiplied according to 
(\ref{multipl}). Thus we have $(x_1 \tensor 1)(1 \tensor 
x_2)(x_3 \tensor 1)(1 \tensor x_4)(x_5 \tensor 1) = 
(\nabla^3 \tensor \nabla^2) \phi (x_1 \tensor x_2 \tensor 
x_3 \tensor x_4 \tensor x_5)$, where $\phi = \tau_3 \tau_4 
\tau_2$ as defined by the given braid diagram.

 We prove now by induction on $n$ that for every product 
$(x_1 \tensor 1)(1 \tensor x_2)(x_3 \tensor 1)\ldots $ with 
$i$ factors of the form $x_j \tensor 1$ and $n-i$ factors of 
the form $1 \tensor x_j$ there is an element $\phi \in B_n$ 
such that 
 $$(x_1 \tensor 1)(1 \tensor x_2)(x_3 \tensor 1)\ldots  = 
(\nabla^i \tensor \nabla^{n-i}) \phi(x_1 \tensor \ldots 
\tensor x_n).$$
 Furthermore if $t$ denotes the number of pairs of factors 
$f_1,f_2$ in the product $(x_1 \tensor 1)(1 \tensor x_2)(x_3 
\tensor 1)\ldots $ where $f_1$ is to the left of $f_2$, 
$f_1$ is of the form $(1 \tensor x_j)$ and $f_2$ is of the 
form $(x_j \tensor 1)$, or briefly the number of factors in 
reverse position, then $\phi$ is composed of $t$ generators 
$\tau_j$ of $B_n$. Observe that $\phi$ and the number $t$ 
are uniquely determined by the properties of the 
multiplication of $A \tensor A$ and the braid group $B_n$, 
which has homogeneous relations. 

 For $n = 1$ we have the trivial cases $x \tensor 1 = 
(\nabla^1 \tensor \nabla^0)(x)$ and $1 \tensor x = (\nabla^0 
\tensor \nabla^1)(x)$, where $\nabla^1 = \id$ and $\nabla^0 
= 1$. For the induction nothing is to be proved if $i = n$ 
or $i = 0$. In these cases we have $t = 0$.

 We assume now that the claim is true for $n$. The induction 
step for $i \not= 0,n+1$ is given by 
 $$\begin{array}{r @{\ } l}
 (x_1 \tensor 1)&\cdot (1 \tensor x_2) \cdot (x_3 \tensor 1) 
\cdot \ldots  \cdot (1 \tensor x_{n+1}) =\\
 &= \{(\nabla^i \tensor \nabla^{n-i})\phi(x_1 \tensor \ldots 
\tensor x_n)\}\cdot (1 \tensor x_{n+1}) \\
 &= \{(\nabla^i \tensor \nabla^{n-i}) \sum_k (u_{k,1} \tensor 
\ldots \tensor u_{k,n})\}\cdot (1 \tensor x_{n+1}) \\
 &= (\sum_k u_{k,1} \cdot \ldots \cdot u_{k,i} \tensor 
u_{k,i+1} \cdot \ldots \cdot u_{k,n})
\cdot (1 \tensor x_{n+1}) \\
 &= \sum_k (u_{k,1} \tensor 1) \cdot \ldots \cdot (u_{k,i} 
\tensor 1) \cdot (1 \tensor u_{k,i+1}) \cdot \ldots \cdot 
(1 \tensor u_{k,n}) \cdot (1 \tensor x_{n+1}) \\
 &= (\nabla^i \tensor \nabla^{n-i+1}) \sum_k (u_{k,1} \tensor 
\ldots \tensor u_{k,n} \tensor x_{n+1}) \\
 &= (\nabla^i \tensor \nabla^{n-i+1}) (\phi \tensor 1) (x_1 \tensor 
\ldots \tensor x_{n+1})  
 \end{array} $$ 
 where $t$, the number of factors in reverse position, does 
not change, neither does the number of generators $\tau_i$ 
used in the representation of $\phi \tensor 1$. The second 
possibility is 
 $$\begin{array}{r @{\ } l}
 (x_1 \tensor 1)&\cdot (1 \tensor x_2) \cdot (x_3 \tensor 1) 
\cdot \ldots  \cdot (x_{n+1} \tensor 1) \\
 &= \{(\nabla^i \tensor \nabla^{n-i})\phi(x_1 \tensor \ldots 
\tensor x_n)\}\cdot (x_{n+1} \tensor 1) \\
 &= (\sum_k u_{k,1} \cdot \ldots \cdot u_{k,i} \tensor 
u_{k,i+1} \cdot \ldots \cdot u_{k,n})
\cdot (x_{n+1} \tensor 1) \\
 &= (\sum_k u_{k,1} \cdot \ldots \cdot u_{k,i} \tensor 
u_{k,i+1} \cdot \ldots \cdot u_{k,n-1}) \cdot (1 \tensor 
u_{k,n}) \cdot (x_{n+1} \tensor 1) \\
 \end{array} $$ 

 $$\begin{array}{r @{\ } l}
 &= (\sum_k u_{k,1} \cdot \ldots \cdot u_{k,i} \tensor 
u_{k,i+1} \cdot \ldots \cdot u_{k,n-1}) \cdot (v_{k,n} \tensor 
1) \cdot (1 \tensor v_{k,n+1}) \\
 &= (\nabla^{i+1} \tensor \nabla^{n-i-1})\rho(\sum_k u_{k,1} 
\tensor \ldots \tensor u_{k,n-1} \tensor v_{k,n}) \cdot (1 
\tensor v_{k,n+1}) \\ 
 &= (\nabla^{i+1} \tensor \nabla^{n-i})(\rho \tensor 
1)(\sum_k u_{k,1} \tensor \ldots \tensor u_{k,n-1} \tensor 
v_{k,n} \tensor v_{k,n+1}) \\ 
 &= (\nabla^{i+1} \tensor \nabla^{n-i})(\rho \tensor 1)(1^{n-
1} \tensor \tau)(\phi \tensor 1)(x_1 \tensor \ldots 
\tensor x_{n+1}).
 \end{array} $$ 
 where $\phi(x_1 \tensor \ldots \tensor x_n) = \sum_k 
u_{k,1} \tensor \ldots \tensor u_{k,n}$, $\tau( u_{k,n} 
\tensor x_{n+1}) = \sum v_{k,n} \tensor v_{k,n+1}$, and 
$(1^{n-1} \tensor \tau)(\phi \tensor 1)(x_1 \tensor \ldots 
\tensor x_n \tensor x_{n+1}) = \sum_k u_{k,1} \tensor \ldots 
\tensor u_{k,n-1} \tensor v_{k,n} \tensor v_{k,n+1}$. We 
determine the number $t(\psi)$ of generators $\tau_i$ 
occurring in $\psi = (\rho \tensor 1)(1^{n-1} \tensor 
\tau)(\phi \tensor 1)$. We have by induction $t(\phi) = t_n$ 
the number of factors in $(x_1 \tensor 1)\cdot (1 \tensor 
x_2) \cdot (x_3 \tensor 1) \cdot \ldots $ in reverse 
position. Also we have $t_{n+1} = t_n + (n-i)$ the number of 
factors in $(x_1 \tensor 1)\cdot (1 \tensor x_2) \cdot (x_3 
\tensor 1) \cdot \ldots \cdot (x_{n+1} \tensor 1)$ in 
reverse position. Then $t(\psi) = t((\rho \tensor 1)(1^{n-1} 
\tensor \tau)(\phi \tensor 1)) = t(\rho \tensor 1) + t(1^{n-
1} \tensor \tau) + t(\phi \tensor 1) = (n-i-1) + 1 + t_n = 
t_{n+1}$.

 If we sum up we obtain
 $$(x_1 \tensor 1 + 1 \tensor x_1) \cdot \ldots \cdot (x_n 
\tensor 1 + 1 \tensor x_n) = \sum_i \sum_{\phi_i} (\nabla^i 
\tensor \nabla^{n-i}) \phi_i(x_1 \tensor \ldots \tensor 
x_n),$$
 for certain $\phi_i \in B_n$ which arise in the evaluation 
given above. 

 Now let $z \in A^n(\zeta)$. We expand the products in 
$\nabla^n p^n z = \sum_k (x_{k,1} \tensor 1 + 1 \tensor 
x_{k,1}) \cdot \ldots \cdot (x_{k,n} \tensor 1 + 1 \tensor 
x_{k,n})$. Each of these products in the sum is treated in 
the same way as described above. Using (\ref{zetasymm}) we get 
 $$\begin{array}{rl}
 \nabla^n p^n(z) &= \sum_k (x_{k,1} \tensor 1 + 1 \tensor 
x_{k,1}) \cdot \ldots \cdot (x_{k,n} \tensor 1 + 1 \tensor 
x_{k,n})\\
 &= \sum_k \sum_i \sum_{\phi_i} (\nabla^i \tensor \nabla^{n-
i}) \phi_i (x_{k,1} \tensor \ldots \tensor x_{k,n}) \\
 &= \sum_i \sum_{\phi_i} (\nabla^i \tensor \nabla^{n-
i}) \phi_i (z) \\
 &= \sum_i \sum_{\phi_i} (\nabla^i \tensor \nabla^{n-
i}) \zeta^{-t(\phi_i)} \rho_i (z) \\
 \end{array}$$
 where $\rho_i \in S_n$ are the canonical images of the 
$\phi_i \in B_n$ and $t(\phi_i)$ is the number of factors 
$\tau_j$ in the representation of $\phi_i$.

 This gives us 
 $$\begin{array}{rl}
 [p^n(z)] &= \sum_{\sigma \in S_n} \nabla^n p^n \sigma(z) \\
 &= \sum_\sigma \sum_i \sum_{\phi_i} \zeta^{-t(\phi_i)} 
(\nabla^i \tensor \nabla^{n-i}) \rho_i\sigma (z)\\ 
 &= \sum_\sigma \sum_i \big(\sum_{\phi_i} \zeta^{-t(\phi_i)} 
\big) (\nabla^i \tensor \nabla^{n-i}) \sigma (z) \\
 &= \sum_i c_i (\nabla^i \tensor \nabla^{n-i}) \sum_\sigma 
\sigma(z)
 \end{array}$$
 where the factors $c_i = \sum_{\phi_i} \zeta^{-t(\phi_i)} 
\in k$. We want to show that the $c_i$ are zero for all $0 
< i < n$. 

 So fix $n$ and $i$. Consider one product $(x_1 \tensor 
1)\cdot (1 \tensor x_2) \cdot (x_3 \tensor 1) \cdot \ldots 
$ in the development of $(x_1 \tensor 1 + 1 \tensor x_1) 
\cdot \ldots \cdot (x_n \tensor 1 + 1 \tensor x_n) = \sum_i 
\sum_{\phi_i} (\nabla^i \tensor \nabla^{n-i}) \phi_i(x_1 
\tensor \ldots \tensor x_n)$ and its corresponding 
$\phi_i$. The chosen summand is completely determined by 
giving the positions in $\{1, \ldots, n\}$ of the $n-i$ 
factors of the form $(1 \tensor x_j)$. The first of these 
factors has $\lambda_1$ factors of the form $(x_j \tensor 
1)$ to its right with $0 \leq \lambda_1 \leq i$. So it 
contributes $\lambda_1$ pairs of factors in reverse 
position. The second factor of the form $(1 \tensor x_j)$ 
contributes $\lambda_2$ (with $0 \leq \lambda_2 \leq 
\lambda_1 \leq i$) pairs of factors in reverse position, 
and so on. We obtain $t = \lambda_1 + \lambda_2 + \ldots + 
\lambda_{n-i}$ pairs in reverse position. If we know the 
$\lambda_i$ with $0 \leq \lambda_{n-i} \leq \ldots \leq 
\lambda_2 \leq \lambda_1 \leq i$ then they also determine 
uniquely the position of the factors of the form $(1 
\tensor x_j)$. Each partition of $t = \lambda_1 + \lambda_2 
+ \ldots + \lambda_{n-i}$  into (at most) $n-i$ parts each 
$\leq i$ gives one term $\zeta^{-t}$ in $c_i = 
\sum_{\phi_i} \zeta^{-t(\phi_i)}$ and we find $p(i,n-i,t)$ 
partitions of $t$ into at most $n-i$ parts each $\leq i$. 
So we get 
 $$c_i = \sum_{t \geq 0} p(i,n-i,t) \zeta^{-t}.$$ 
 By a theorem of Sylvester (\cite {A} Theorem 3.1) we have 
 $$\sum_{t \geq 0} p(i,n-i,t)q^t = \frac{(1-q^n)(1-q^{n-1}) 
\ldots (1-q^{n-i+1})}{(1-q^i)(1-q^{i-1}) \ldots (1-q)}$$ 
 hence $c_i = 0$ for $0 < i < n$ since $\zeta$ and also 
$\zeta^{-1}$ are {\em primitive} $n$-th roots of unity. 

 So we have shown 
 $$\begin{array}{rl}
 [p^n(z)] &= \sum_{\sigma \in S_n} \nabla^n p^n \sigma(z) \\
 &= \sum_i c_i (\nabla^i \tensor \nabla^{n-i}) \sum_\sigma 
\sigma(z)\\
 &= \sum_\sigma \nabla^n \sigma(z) \tensor 1 + 1 \tensor 
\sum_\sigma \nabla^n \sigma(z) \\
 &= p[z].  
 \end{array}$$
 \end{pf}

 \begin{cor}
 Let $H$ be a Hopf algebra in $\C$. Then the set of 
primitive elements $P(H)$ forms a Lie algebra in $\C$. 
 \end{cor}

 \begin{pf}
 By Lemma \ref{primYD} $P(H)$ is a Yetter-Drinfeld 
submodule of $H$. Let $z \in P(H)^n(\zeta)$. Then $p([z]) = 
[p^n(z)] = [\Delta^n(z)] = \Delta([z])$ since $\Delta$ is 
an algebra homomorphism. Hence $[z] \in P(H)$. So $P(H)$ is 
a Lie subalgebra of $H^L$. 
 \end{pf}

 \begin{defn}{\em
 Let $A$ be an algebra in $\C$ and let $\rend(A)$ be the 
inner endomorphism object of $A$ in $\C$, i.e. the Yetter-
Drinfeld module $\rend(A)$ satisfying $\C(X \tensor A, A) 
\iso \C(X, \rend(A))$ for all $X \in \C$. It can be shown 
that 
 $$\begin{array}{r@{}l}
 \rend(A) := \{f \in \Hom(A,A) \vert &\exists \sum f_{(0)} 
\tensor f_{(1)} \in \Hom(A,A) \tensor K \forall a \in A:\\ 
&\sum f_{(0)}(a) \tensor f_{(1)} = \sum f(a_{(0)})_{(0)} 
\tensor f(a_{(0)})_{(1)} S(a_{(1)}) \} 
 \end{array}$$
 is the Yetter-Drinfeld module with the required universal 
property. $\rend(A)$ operates on $A$ by a canonical map 
$\ev: \rend(A) \tensor A \to A$ with $\ev(f \tensor a) = 
f(a)$. 

 A {\it derivation} from $A$ to $A$  is a linear map
$(d: A \to A) \in \rend(A)$ such that
 $$d(ab) = d(a)b + (1 \tensor d)(\tau \tensor 1)(d \tensor 
a \tensor b)$$
 for all $a,b \in A$. Observe that in the symmetric 
situation this means $d(ab) = d(a)b + ad(b)$.} 
 \end{defn}

 It is clear that all derivations from $A$ to $A$ form an 
object $\Der(A)$ in $\C$ and that there is an operation 
$\Der(A) \tensor A \longrightarrow A$. 

 \begin{cor}
 $\Der(A)$ is a Lie algebra.
 \end{cor}

 \begin{pf} 
 Let $m$ denote the multiplication of $A$. An endomorphism 
$x: A \longrightarrow A$ in $\rend(A)$ is a derivation iff 
$m(x \tensor 1 + 1 \tensor x) = xm$ where $(x \tensor y)(a 
\tensor b) = (\ev \tensor \ev)(1 \tensor \tau \tensor 1)(x 
\tensor y \tensor a \tensor b)$ for elements $a$ and $b$ in 
$A$ and elements $x$ and $y$ in $\rend(A)$. So $x \in 
\rend(A)$ is a derivation iff $m p(x) = xm$. 

 To show that $\Der(A)$ is a Lie algebra it suffices to 
show that it is closed under Lie multiplication since it is 
a subobject of $\rend(A)$, which is an algebra in the 
category $\C$. Let $\zeta$ be a primitive $n$-th root of 
unity. Let $\nabla: \rend(A) \tensor \rend(A) \to \rend(A)$ 
be the multiplication of $\rend(A)$. 

 If $x_1,x_2 \in \Der(A)$ then $m p(x_1)p(x_2) = x_1mp(x_2) 
= x_1x_2m$ or more generally $m (\nabla^n p^n)(x_1 \tensor 
\ldots \tensor x_n) = \nabla^n (x_1 \tensor \ldots \tensor 
x_n)m$ for all $x_1 \tensor \ldots \tensor x_n \in 
\Der(A)^n$. Thus we get for $z \in \Der(A)^n(\zeta)$ 
 $$\begin{array}{r@{}l}
 mp([z]) &= m [p^n(z)] = \sum m (\nabla^n \sigma (p^n 
(z)))\\ 
 &= \sum m (\nabla^n p^n \sigma(z)) = \sum \nabla^n \sigma(z) 
m = [z]m 
 \end{array}$$
 hence $[z] \in \Der(A)$.
 \end{pf}

 \section{The Universal Enveloping Algebra of a Lie Algebra}

 As in \cite{P1} we can now construct the universal 
enveloping algebra of a Lie algebra $P$ in $\cal C$ as 
$U(P) := T(P)/I$ where $T(P)$ is the tensor algebra over 
$P$, which lives again in $\cal C$, and where $I$ is the 
ideal generated by the relations $[z] - \sum \nabla^n 
\sigma(z) $ for all $z \in P^n(\zeta)$, for all $n$ and for 
all primitive $n$-th roots of unity $\zeta$. Then $U(P)$ 
clearly is a universal solution for the following universal 
problem 
  $$\bfig
 \putmorphism(0, 400)(1, 0)[P`U(P)` \iota]{800}1a
 \putmorphism(0, 400)(2, -1)[`` f]{800}1l
 \putmorphism(800, 400)(0, 1)[`A` g]{400}1r
 \efig$$
 where for each morphism of Lie-algebras $f$ there is a 
unique morphism of algebras $g$ such that the diagram 
commutes. 

 \begin{thm}
 Let $P$ be a Lie algebra in $\C$. Then the universal 
enveloping algebra $U(P)$ is a Hopf algebra in $\C$. 
 \end{thm}

 \begin{pf}
 It is easily seen that $\delta: P \to (U(P) \tensor 
U(P))^L$ in $\M^{kG}$ given by $\delta(x) := \overline x 
\tensor 1 + 1 \tensor \overline x$ where $\overline x$ is 
the canonical image of $x \in P$ in $U(P)$ and the counit 
$\varepsilon: U(P) \to k$ given by the zero morphism $0: P 
\to k$ define the structure of a bialgebra on $U(P)$ in 
$\C$. 

 Now we want to define $S: U(P) \to U(P)^{op+}$ by the Lie 
homomorphism $S: P \to U(P)^{op+}$, $S(x) = -\bar x$. Here 
$A^{op+}$ is the algebra obtained from the algebra $A$ by 
the multiplication $A \tensor A \buildrel \tau \over 
\longrightarrow A \tensor A \buildrel \nabla \over 
\longrightarrow A$. Then for $z \in P^n(\zeta)$ we have 
 $$ \begin{array}{rl}
 S([z])
 &= -\overline{[z]}
 = -\sum_\sigma \nabla^n\sigma(\overline z)
 = -\sum_\sigma \nabla^n \pi {\pi}^{-1} \sigma(\overline z)\\
 &= -\sum_\sigma (\nabla^{op})^n {\pi}^{-1} 
\sigma(\overline z) \mbox{ (by (3)) } 
 = -\sum_\sigma (\nabla^{op})^n \zeta^{\frac{n(n-1)}2} {\rho}^{-1} 
\sigma(\overline z)\\
 &= -\zeta^{\frac{n(n-1)}2}[\overline z]
 = (-1)^n[\overline z]
 = [\overline{S^n(z)}]
 \end{array}$$ 
 where $\pi \in B_n$ is the braid map given by the twist of 
all $n$ strands with source $\{1,\ldots,n\}$ and domain 
$\{n, \ldots,1\}$, $\pi = (\tau_1) \ldots (\tau_{n-2} 
\ldots \tau_2 \tau_1) (\tau_{n-1} \ldots \tau_2 \tau_1)$ 
and\break $\zeta^{- \frac{n(n-1)}2}\pi = \rho \in S_n$. 

 Hence $S$ is a Lie homomorphism and factorizes through 
$U(P)$. Since $U(P)$ is generated as an algebra by $P$ we 
prove that $S$ is the antipode by complete induction: 
 $$\nabla(1 \tensor S)\Delta(1) = 1S(1) = 1 = \varepsilon (1), $$ 
 $$\nabla(1 \tensor S)\Delta(\overline x) = \overline x + 
S(\overline x) = 0 = \varepsilon(\overline x).$$ 

 Before we prove the general induction step we observe that 
$\Delta: U(P) \to U(P) \tensor U(P)$ is a morphism in $\C = 
\YD_K^K$ so that we have in particular 
 $$\sum (a^0)_1 \tensor (a^0)_2 \tensor a^1 = \sum (a_1)^0 
\tensor (a_2)^0 \tensor (a_1)^1(a_2)^1 \in U(P) \tensor 
U(P) \tensor K$$ 
 for $a \in U(P)$. (Here we use $\delta(a) = \sum a^0 
\tensor a^1$ to denote the comodule structure in 
$\YD^K_K$.) Assume now that $a$ is writte as a product of 
$n \geq 1$ elements in $P$ and that $\sum a_1S(a_2) = 0$. 
Then for all $x \in P$ we have $\sum 
(a_1)^0S((a_2)^0)S(x(a_1)^1(a_2)^1) = \sum 
(a^0)_1S((a^0)_2)S(xa^1) = 0$ since $\delta(a) = \sum a^0 
\tensor a^1 \in \overline{P \tensor \ldots \tensor P} 
\tensor K \subseteq U(P) \tensor K$. So we have 
 $$\begin{array}{rl}
 \nabla (1 \tensor S) \Delta (xa)
 &= \nabla (1 \tensor S) \sum (xa_1 \tensor a_2 + (a_1)^0 
\tensor (x(a_1)^1)a_2) \\ 
 &= \sum xa_1S(a_2) + \sum (a_1)^0S((x(a_1)^1)a_2) \\
 &= \sum (a_1)^0 S((a_2)^0) S(x(a_1)^1 (a_2)^1) = 0 = \eta 
\varepsilon (xa). 
 \end{array}$$
 The second condition $\nabla (S \tensor 1) \Delta = \eta 
\varepsilon$ is proved in a similar way (by using elements 
of the form $ax$ and the equation $\sum S((a \kappa)_1)(a 
\kappa)_2 = 0$ for $a$ written as a product of $n$ elements 
in $P$ and $\kappa \in K$). So $S$ is an antipode and 
$U(P)$ is a Hopf algebra in $\C$. 
 \end{pf}

 \section{$(G,\chi)$ Lie algebras}

 In \cite{P1} we introduced and studied the concept of 
 $G$-graded Lie algebras or $(G,\chi)$-Lie algebras for an 
abelian group $G$ with a bicharacter $\chi$ generalizing 
the concepts of Lie algebras, Lie super algebras, and Lie 
color algebras. The reader may find examples of such 
$(G,\chi)$-Lie algebras in \cite{P1}. A generalization of 
this concept of Lie algebras to the group graded case for a 
noncommutative group requires the use of Yetter-Drinfeld 
modules over $kG$. We show that  $(G,\chi)$-Lie algebras 
are Lie algebras on Yetter-Drinfeld modules in the sense of 
this paper. We use the notation of \cite{P1}. 

 Let $G$ be an abelian group with a bicharacter $\chi: G 
\tensor_{\Bbb Z} G \to k^*$. Let $P$ be a $kG$-comodule. Then 
$P$ is a Yetter-Drinfeld module over $kG$ \cite{M} with 
the module structure $x \cdot g = \chi(h,g) x$ for 
homogeneous elements $x = x_h \in M$ with $\delta(x) = x \tensor 
h$. The braid map is $\tau(x_h \tensor y_g) = y_g \tensor 
x_h \cdot g = \chi(h,g) y_g \tensor x_h$, hence the 
braiding given in \cite{P1} after Example 2.3.

  Let $\zeta \in k^*$ be given. Let $(g_1, \ldots, g_n)$ be 
a $\zeta$-family, i.e. $\chi(g_i,g_j)\chi(g_j,g_i) = 
\zeta^2$. Let $Q := \sum_{\sigma \in S_n} P_{g_{\sigma(1)}} 
\tensor \ldots \tensor P_{g_{\sigma(n)}}$. 

 \begin{lma}\label{Sn-action}
 $Q$ is a right $S_n$-module by
 $$(x_1 \tensor \ldots \tensor x_n)\sigma = \rho(\sigma, 
(g_1, \ldots, g_n)) x_{\sigma(1)} \tensor \ldots \tensor 
x_{\sigma(n)}$$
 for $\sigma \in S_n$ and $x_1 \tensor \ldots \tensor x_n 
\in P_{g_1} \tensor \ldots \tensor P_{g_n}$.
 \end{lma}

 \begin{pf}
 We have to show the compatibility of this operation with 
the composition of permutations. Let $\sigma, \tau \in 
S_n$. We use Lemma 2.2 of \cite{P1}. Then 
 $$\begin{array}{r@{}l}
 (x_1 \tensor &\ldots \tensor x_n)(\sigma\tau) = \\
 &= \rho(\sigma\tau, (g_1, \ldots, g_n)) x_{\sigma\tau(1)} 
\tensor \ldots \tensor x_{\sigma\tau(n)}\\
 &= \rho(\sigma, (g_1, \ldots, g_n)) \rho(\tau, 
(g_{\sigma(1)}, \ldots, g_{\sigma(n)})) x_{\sigma\tau(1)} 
\tensor \ldots \tensor x_{\sigma\tau(n)}\\
 &= (\rho(\sigma, (g_1, \ldots, g_n))  x_{\sigma(1)} 
\tensor \ldots \tensor x_{\sigma(n)})\tau\\
 &= ((x_1 \tensor \ldots \tensor x_n)\sigma)\tau.  \\
 \end{array}$$
 \end{pf}

 $Q$ becomes a left $S_n$-module by $\sigma(x_1 \tensor 
\ldots \tensor x_n) = \rho(\sigma^{-1}, (g_1, \ldots, g_n)) 
x_{\sigma^{-1}(1)} \tensor \ldots \tensor x_{\sigma^{-
1}(n)}$. Thus $\bigoplus_{\{(g_1,\ldots,g_n) 
\ \zeta\mbox{\small -family}\}} P_{g_1} \tensor 
\ldots \tensor P_{g_n}$ is also a left $S_n$-module. 

 This action is connected with the action of $B_n$ on 
$\bigoplus_{\{(g_1,\ldots,g_n) \ \zeta\mbox{\small -
family}\}} P_{g_1} \tensor \ldots \tensor P_{g_n}$ by
 \begin{equation}\label{Bn-action}
 \zeta^{-1}\tau_i(x_1 \tensor \ldots \tensor x_n) = 
\sigma_i(x_1 \tensor \ldots \tensor x_n)
 \end{equation} 
 for the canonical generators $\tau_i$ of $B_n$ resp. 
$\sigma_i$ of $S_n$, since
  $$\begin{array}{r@{}l}
 \zeta^{-1}\tau_i&(x_1 \tensor \ldots \tensor x_n) =\\
 &= \zeta^{-1}\chi(g_i,g_{i+1}) x_1 \tensor \ldots \tensor 
x_{i+1} \tensor x_i \tensor \ldots \tensor x_n\\
 &= \rho(\sigma_i^{-1},(g_1, \ldots, g_n)) x_{\sigma_i^{-1}(1)} 
\tensor \ldots \tensor x_{\sigma_i^{-1}(n)} \\
 &= \sigma_i(x_1 \tensor \ldots \tensor x_n).
 \end{array} $$
 In particular we have
  $$\begin{array}{r@{}l}
 \tau_i^{-1} &\tau_{i+1}^{-1} \ldots \tau_{j-1}^{-1} 
\tau_{j}^2 \tau_{j-1} \ldots \tau_{i+1} \tau_i(x_1 \tensor 
\ldots \tensor x_n) = \\ 
 &= \zeta^2  \sigma_i^{-1} \sigma_{i+1}^{-1} 
\ldots \sigma_{j-1}^{- 1} \sigma_{j}^2 \sigma_{j-1} \ldots 
\sigma_{i+1} \sigma_i(x_1 \tensor \ldots \tensor x_n)\\
 &= \zeta^2 (x_1 \tensor \ldots \tensor x_n),
 \end{array}$$
 so that $x_1 \tensor \ldots \tensor x_n \in P^n(\zeta)$ 
by Lemma \ref{append1}. Thus we have 
 $$\bigoplus_{\{(g_1,\ldots,g_n) 
\ \zeta\mbox{\small -family}\}} P_{g_1} \tensor 
\ldots \tensor P_{g_n} \subseteq P^n(\zeta).$$
 Conversely let $\sum x_1 \tensor \ldots \tensor x_n \in 
P^n = \bigoplus _{\{(g_1, \ldots, g_n)\}} 
P_{g_1} \tensor \ldots \tensor P_{g_n}$ with homogeneous 
summands and assume that one of the summands is non-zero in 
$ P_{g_1} \tensor \ldots \tensor P_{g_n}$ where $(g_1, 
\ldots, g_n)$ is not a $\zeta$-family for example by 
$\chi(g_i,g_{i+1}) \chi(g_{i+1},g_i) \not= \zeta^2$. Then 
$(\tau_i^2 - \zeta^2)(\sum x_1 \tensor \ldots \tensor x_n)$ 
has a non-zero component in $P_{g_1} \tensor \ldots \tensor 
P_{g_n}$, hence $\sum x_1 \tensor \ldots \tensor x_n$ 
cannot be in $P^n(\zeta)$. This proves

 \begin{prop}
 Let $\zeta \in k^*$ be given. Then 
 $$P^n(\zeta) = \bigoplus_{\{(g_1,\ldots,g_n) 
\ \zeta\mbox{\small -family}\}} P_{g_1} \tensor 
\ldots \tensor P_{g_n}.$$ 
 \end{prop}

 By Lemma \ref{Sn-action} and (\ref{Bn-action}) the bracket 
multiplication of \cite{P1} is a special case of the 
bracket multiplication of this paper and $(G,\chi)$-Lie 
algebras are Lie algebras over Yetter-Drinfeld modules. 

 \begin{eg}{\em
 As a new example of Lie algebras we give one 
family of examples of $(G,\chi)$-Lie algebras. 
Let $G = C_3 = \{ 0,  1,  2\}$ be the cyclic 
group with 3 elements. Define the structure of a right 
$kG$-module on a $kG$-comodule $V$ (i.e. on a $C_3$-graded 
vector space $V = V_0 \oplus V_1 \oplus V_2$) using the 
bicharacter $\chi: C_3 \tensor_{\Bbb Z} C_3 \iso C_3 \to 
k^*$, $\chi( 1 \tensor  1) = \xi$ a primitive
 3-rd root of unity, by $v \cdot g := \chi(\deg(v) \tensor 
g) v = \chi(\deg(v), g) v$ for $g \in G$ and homogeneous 
elements $v \in V$. Then $V$ is a Yetter-Drinfeld module. 

 Let $A := \rend(V)$ be the inner endomorphism object of 
$V$ in $kG\comod$. By Corollary \ref{algebras_are_Lie} $A$ 
is a Lie algebra. One verifies easily (see \cite{P1}) that 
the only non-zero components $A^n(\zeta)$ for the partial 
Lie multiplication are 
 $$A^2(-1) = A_0 \tensor (A_1 \oplus A_2) \oplus 
(A_0 \oplus A_1 \oplus A_2) \tensor A_0$$
 and 
 $$A^3(\xi) = A_1 \tensor A_1 \tensor A_1 \oplus A_2 
\tensor A_2 \tensor A_2.$$ 

 Now let $\langle .,. \rangle : V \tensor V \to k$ be a 
bilinear form on $V$ in $\C$. We define
 $$\og(V)_i := \{f \in A_i| \forall v,w \in V, \deg(v) = 
j:\ \langle f(v), w \rangle = -\chi( i, 
j) \langle v, f(w) \rangle\}.$$
 This space is the homogeneous component of $\og(V) 
\subseteq A$ that becomes a Yetter-Drinfeld module. 

 For $f \in \og(V)_0$ and $g \in \og(V)_i$, $i \in C_3$, $v 
\in V_j$, $w \in V$ we have
 $$\begin{array}{r@{}l}
 \langle [f,g](v), w \rangle &= \langle (fg - gf)(v), w 
\rangle\\
 &= \langle fg(v), w\rangle - \langle gf(v), w \rangle\\
 &= \chi( i, j)\langle v, gf(w) \rangle - \chi( 
i,  j) \langle v, fg(w) \rangle\\
 &= -\chi( i,  j) \langle v, [f,g](w) \rangle,
 \end{array}$$
hence $[f,g] \in \og(V)_i$. Analogously one shows  $[g,f] 
\in \og(V)_i$.

 For $k = 1,2,3$ let $f_k \in \og(V)_i$ ($i = 1$ or $i = 
2$). Then 
 $$\begin{array}{r@{}l}
 \langle [f_1,f_2,f_3](v), &w \rangle = \sum_{\sigma \in 
S_3} \langle f_{\sigma(1)} f_{\sigma(2)} f_{\sigma(3)} 
(v), w \rangle\\
 &= (-1) \sum_{\sigma \in S_3}
 \chi( i,  i +  i +  j)
 \chi( i,  i +  j)
 \chi( i,  j)
 \langle v, f_{\sigma(3)} f_{\sigma(2)} f_{\sigma(1)} (w) 
\rangle\\ 
 &= -\langle v, [f_1,f_2,f_3](w) \rangle,
 \end{array}$$
 hence $[f_1,f_2,f_3] \in \og(V)_0$. Thus we 
have a Lie algebra $\og(V)$. Depending on the choice of the 
bilinear form this is a generalization of the orthogonal or 
the symplectic Lie algebra.}
 \end{eg}

 \section{Appendix} \label{appendix}

 \begin{pf} of Lemma \ref{append1}:

 Define actions $\pi_{i,j}$ for $1 \leq i < j \leq n$ on $M$ 
by 
 \begin{equation}\label{form2}
 \pi_{i,j} := \tau_i^{- 1} \tau_{i+1}^{- 1} 
\ldots \tau_{j-2}^{- 1} \tau_{j-1} \tau_{j-2} \ldots 
\tau_{i+1} \tau_i\end{equation} 
 Observe that $\pi_{i,i+1} = \tau_i$. 
 Since $ \tau_i \tau_j  = \tau_j \tau_i  \mbox{ if } |i-j| 
\geq 2$ a simple calculation gives
 \begin{equation}\label{form3}
 \begin{array}{c}
 \pi_{i,j} \tau_k = \tau_k \pi_{i,j} \mbox{ for all } k < 
i-1 \mbox{ and all } k > j, \\ 
 \pi_{i,j} \tau_{i-1} = \tau_{i-1} \pi_{i-1,j}, \\
 \pi_{i,j} \tau_k = \tau_k \pi_{i,j} \mbox{ for all } i < k 
< j-1,\\ 
 \pi_{i,j} \tau_{j-1} = \tau_{j-1} \pi_{i,j-1} \mbox{ if } 
i < j-1 \mbox{ and }\\ 
 \pi_{i,j} \tau_{j-1} = \tau_i \tau_i = \tau_{j-1} 
\pi_{i,j} \mbox{ if } i=j-1. 
 \end{array}
 \end{equation} 

 Let $N \subseteq M$ be a $kB_n$ submodule of $M$. Assume 
furthermore that $\tau_i^2 \tau_{i+1} = \tau_{i+1} \tau_i^2$ 
on $N$ for all $i=1, \ldots, n-2$. Then $ \tau_{i+1}^2 
\tau_i = \tau_i^3 = \tau_i \tau_{i+1}^2.$ Consequently we 
have 
 \begin{equation}
 \tau_j^2 \tau_i = \tau_i \tau_j^2
 \end{equation}
 on $N$ for all $i,j = 1, \ldots, n-1.$ Thus the $\tau_j^2$ 
commute with all $\phi \in B_n$ if they act on $N$.

 We introduce the vector subspace $M(\zeta) \subseteq 
\overline{M(\zeta)} \subseteq M$ by 
 $$\overline{M(\zeta)} := \{z \in M | \forall 1 \leq i < j 
\leq n : \pi_{i,j}^2(z) = \zeta^2 z \}$$ 
 and show that $\overline{M(\zeta)}$ is invariant under the 
action of the $\tau_i$ and $\tau_i^2 \tau_{i+1} = 
\tau_{i+1} \tau_i^2$ on $M(\zeta)$ for all $i=1, \ldots, n-
2$. 

 For $z \in M(\zeta)$ and $i < j$ we have $\pi_{i,j}^2 
\tau_k(z) = \tau_k \pi_{i,j}^2(z) = \zeta^2\tau_k(z)$ for 
all $k$ with $1 \leq k < i-1$ and $j < k \leq n$ by 
(\ref{form3}) and for all $k$ with $i < k < j-1$ by 
(\ref{form3}). Furthermore we have $\pi_{i,j}^2 \tau_{i-
1}(z) = \tau_{i-1} \pi_{i-1,j}^2(z) = \zeta^2\tau_{i-1}(z)$ 
by (\ref{form3}), $\pi_{i,j}^2 \tau_{j-1}(z) = \tau_{j-1} 
\pi_{i,j-1}^2(z) = \zeta^2 \tau_{j-1}(z)$ (for $i < j-1$) by 
(\ref{form3}), and $\pi_{i,j}^2 \tau_{j-1}(z) = \tau_{j-1} 
\pi_{i,j}^2(z) = \zeta^2 \tau_{j-1}(z)$ (for $i = j-1$) by 
(\ref{form3}). So there remain two cases to investigate for 
which we use $\pi_{i,j}^2 (z) = \zeta^2 z$ and symmetrically 
$\pi_{i,j}^{-2} (z) = \zeta^{-2} z$ for all $z \in 
M(\zeta)$.
 
 In the first case we get
 $$\begin{array}{rl}
 \pi_{i,j}^2 \tau_i(z)
 &= \tau_i^{-1} \ldots \tau_{j-1}^2 \ldots \tau_i \tau_i (z) 
 = \tau_i^{-1} \ldots \tau_{j-1}^2 \ldots \tau_{i+1} 
(\zeta^2 z)\\
 &= \zeta^2 \tau_i^{-1} \ldots \tau_{j-1}^2 \ldots 
\tau_{i+1} (z) = \zeta^2 \tau_i^{-1} \pi_{i+1,j}^2 (z)\\ 
 &= \zeta^2 \tau_i^{-1}(\zeta^2 z)
 = \zeta^2 \tau_i^{-1} \tau_i^2 (z)
 = \zeta^2 \tau_i (z)\\
 \end{array}$$
 for $i+1<j$ and
 $\pi_{i,i+1}^2 \tau_i(z) = \tau_i^3(z) = \zeta^2 
\tau_i(z)$.

 In the second case we get
 $$\begin{array}{rl}
 \pi_{i,j}^2 \tau_j(z)
 &= \tau_i^{-1} \ldots \tau_{j-1}^2 \ldots \tau_i \tau_j (z) 
 = \zeta^2 \tau_i^{-1} \ldots \tau_{j-1}^2 \ldots \tau_i 
\tau_j \tau_j^{-2} (z) \\
 &= \zeta^2 \tau_j^{-1} \tau_j \tau_i^{-1} \ldots \tau_{j-
1}^2 \ldots \tau_i \tau_j^{-1} (z) 
 = \zeta^2 \tau_j^{-1} \tau_i^{-1} \ldots \tau_j \tau_{j-
1}^2 \tau_j^{-1} \ldots \tau_i (z) \\
 &= \zeta^2 \tau_j^{-1} \tau_i^{-1} \ldots \tau_{j-1}^{-1} 
\tau_j^2 \tau_{j-1} \ldots \tau_i (z) 
 = \zeta^2 \tau_j^{-1} \pi{i,j+1}^2 (z)\\
 &= \zeta^2 \tau_j^{-1} (\zeta^2 z)
 = \zeta^2 \tau_j (z).\\
 \end{array}$$

 Hence we have $\tau_i(z) \in M(\zeta)$ for all $z \in 
M(\zeta)$ and all $i = 1, \ldots, n-1$. 

 The claim $\tau_i^2 \tau_{i+1} = \tau_{i+1} \tau_i^2$ is 
clear from the invariance and the fact, that $\tau_i^2$ on 
$M(\zeta)$ is multiplication by $\zeta^2$.

 Since the $\tau_i^2$ commute in their action on 
$\overline{M(\zeta)}$ with all $\phi \in B_n$ it is clear 
that $\overline{M(\zeta)} \subseteq M(\zeta)$. 
 \end{pf}

 We now study specific braids. The following identity
 \bgr{110}{80} 
 \put(0,70){\braid} 
 \put(20,50){\braid} 
 \put(30,40){\braid} 
 \put(30,30){\braid} 
 \put(70,40){\braid} 
 \put(70,30){\braid} 
 \put(80,20){\braid} 
 \put(100,0){\braid} 
 \put(0,0){\ibraid} 
 \put(20,20){\ibraid} 
 \put(80,50){\ibraid} 
 \put(100,70){\ibraid} 
 \put(0,10){\idgr{60}} 
 \put(20,30){\idgr{20}} 
 \put(30,0){\idgr{20}} 
 \put(30,60){\idgr{20}} 
 \put(40,0){\idgr{30}} 
 \put(40,50){\idgr{30}} 
 \put(70,0){\idgr{30}} 
 \put(70,50){\idgr{30}} 
 \put(80,0){\idgr{20}} 
 \put(80,60){\idgr{20}} 
 \put(90,30){\idgr{20}} 
 \put(110,10){\idgr{60}} 
 \put(21,0){\objo{...}} 
 \put(91,0){\objo{...}} 
 \put(21,76){\objo{...}} 
 \put(91,76){\objo{...}} 
 \put(55,48){\obju{=}} 
 \put(20,60){\twist{-10}{10}} 
 \put(10,10){\twist{10}{10}} 
 \put(100,10){\twist{-10}{10}} 
 \put(90,60){\twist{10}{10}} 
 \egr 
 implies
 \begin{equation}\label{equn1} 
 \tau_1^{-1} \ldots \tau_{i-1}^{-1} \tau_i^2 \tau_{i-1} 
\ldots \tau_1 = \tau_i \ldots \tau_2 \tau_1^2 \tau_2^{-1} 
\ldots \tau_i^{-1} 
 \end{equation} 
 and similarly $\tau_1 \ldots \tau_{i-1} \tau_i^2 \tau_{i-
1}^{-1} \ldots \tau_1^{-1} = \tau_i^{-1} \ldots \tau_2^{-1} 
\tau_1^2 \tau_2 \ldots \tau_i$ for all $i=1, \ldots, n$. 

 Let $B_n \ni \phi \mapsto \tilde\phi \in S_n$ denote the 
canonical epimorphism. 

 For each braid $\phi \in B_n$ there exists a braid 
$\phi_{(i)} \in B_{n+1}$ such that the diagram 
 $$\bfig
 \putmorphism(0, 700)(1, 0)[P \tensor \ldots \tensor P^2 
\tensor \ldots \tensor P`P \tensor \ldots \tensor P`1 
\tensor \ldots \tensor f \tensor \ldots \tensor 1]{2100}1a
 \putmorphism(0, 600)(1, 0)[\phantom{P \tensor \ldots 
\tensor\ } i \phantom{2 \tensor \ldots \tensor P}`\phantom{P 
\tensor \ldots \tensor P}`]{2100}0a 
 \putmorphism(0, 100)(1, 0)[P \tensor \ldots \tensor P^2 
\tensor \ldots \tensor P`P \tensor \ldots \tensor P`1 
\tensor \ldots \tensor f \tensor \ldots \tensor 1]{2100}1a
 \putmorphism(0, 0)(1, 0)[\phantom{P \tensor \ldots 
\tensor\ } j \phantom{2 \tensor \ldots \tensor P}`\phantom{P 
\tensor \ldots \tensor P}`]{2100}0a 
 \putmorphism(0, 700)(0, 1)[``\phi_{(i)}]{600}1r
 \putmorphism(2100, 700)(0, 1)[``\phi]{600}1r
 \efig$$
 commutes for all $f: P^2 \to P$ in $\C$ (where $j = 
\tilde\phi(i)$). The braid $\phi_{(i)}$ can be given 
explicitly, but we are only interested in the following 
special forms 
 $$\begin{array}{r@{\ }ll @{\qquad} r@{\ }l}
 \tau_{j(i)} &= \tau_{j+1} & \mbox{ if } j > i; & 
 \tau_{i-1(i)} &= \tau_i\tau_{i-1};\\
 \tau_{j(i)} &= \tau_j & \mbox{ if } j < i-1; &
 \tau_{i(i)} &= \tau_i \tau_{i+1} \\ 
 \end{array}$$
 which can be easily verified.

 By (\ref{equn1}) we have for all $z \in P^{n+1}(-1,\zeta)$ 
 \begin{equation}\label{equn2} 
 \tau_i^2 \tau_{i-1} \ldots \tau_1(z) = \tau_{i-1} \ldots \tau_1(z).
 \end{equation}

 \begin{lma}\label{equn5}
 For $z \in P^{n+1}(-1,\zeta)$, $\phi 
\in B_n$ and $j := \tilde{\phi}(i)$ we have 
 $$ \begin{array}{r @{\ } l } 
 \phi_{(i)} \tau_{i-1} \ldots \tau_1(z) &= \tau_{j-1} 
\ldots \tau_1 (1 \tensor \phi)(z) ; \\
 \phi_{(i)} \tau_i \ldots \tau_1(z) &= \tau_j 
\ldots \tau_1 (1 \tensor \phi)(z) . \\
 \end{array}$$
 \end{lma}

 \begin{pf}
 To prove this we first observe that these two relations 
are compatible with the group structure of $B_n$. 
For $\tilde{\phi} \tilde{\psi}(i) = \tilde{\phi}(j) = k$ we 
have
 $$ \begin{array}{r @{\ } l }
 \phi_{(j)} \psi_{(i)} \tau_{i-1} \ldots \tau_1(z)
 &= \phi_{(j)} \tau_{j-1} \ldots \tau_1 (1 \tensor \psi)(z)
 = \tau_{k-1} \ldots \tau_1 (1 \tensor \phi\psi)(z) ; \\ 
 \phi_{(j)} \psi_{(i)} \tau_i \ldots \tau_1(z)
 &= \phi_{(j)} \tau_j \ldots \tau_1 (1 \tensor \psi)(z)
 = \tau_k \ldots \tau_1 (1 \tensor \phi\psi)(z)  \\ 
 \end{array}$$
 so we only have to show these relations for the generators 
$\phi = \tau_j$, $j = 1, \ldots, n-1$. In these cases we 
have 
 $$\begin{array}{r @{\ }ll}
 \tau_{j(i)} \tau_{i-1} \ldots \tau_1 (z)
   &= \tau_{j+1} \tau_{i-1} \ldots \tau_1 (z) &\\
   &= \tau_{i-1} \ldots \tau_1 \tau_{j+1} (z) &\\
   &= \tau_{i-1} \ldots \tau_1 (1 \tensor \tau_j) (z)
   & \mbox{ for } j > i;\\
 \tau_{j(i)} \tau_{i-1} \ldots \tau_1 (z)
   &= \tau_j \tau_{i-1} \ldots \tau_1 (z) &\\
   &= \tau_{i-1} \ldots \tau_j \tau_{j+1} \tau_j \ldots 
    \tau_1 (z) &\\ 
   &= \tau_{i-1} \ldots \tau_{j+1} \tau_j \tau_{j+1} \ldots 
    \tau_1 (z) &\\ 
   &= \tau_{i-1} \ldots \tau_1 \tau_{j+1} (z) &\\
   &= \tau_{i-1} \ldots \tau_1 (1 \tensor \tau_j) (z)
   & \mbox{ for } j < i-1;\\
 \tau_{i-1(i)} \tau_{i-1} \ldots \tau_1 (z)
   &= \tau_i \tau_{i-1} \tau_{i-1} \ldots \tau_1 (z) &\\
   &= \tau_i \tau_{i-2} \ldots \tau_1 (z) &\\
   &= \tau_{i-2} \ldots \tau_1 \tau_i (z) &\\
   &= \tau_{i-2} \ldots \tau_1 (1 \tensor \tau_{i-1}) (z);
   &\\
 \tau_{i(i)} \tau_{i-1} \ldots \tau_1 (z)
   &= \tau_i \tau_{i+1} \tau_{i-1} \ldots \tau_1 (z) &\\
   &= \tau_i \ldots \tau_1 \tau_{i+1}(z) &\\
   &= \tau_i \ldots \tau_1 (1 \tensor \tau_i) (z);
   &\\
 \end{array}$$
 $$\begin{array}{r @{\ }ll}
 \tau_{j(i)} \tau_i \ldots \tau_1 (z)
   &= \tau_{j+1} \tau_i \ldots \tau_1 (z) &\\
   &= \tau_i \ldots \tau_1 \tau_{j+1} (z) &\\
   &= \tau_i \ldots \tau_1 (1 \tensor \tau_j) (z)
   & \mbox{ for } j > i;\\
 \tau_{j(i)} \tau_i \ldots \tau_1 (z)
   &= \tau_j \tau_i \ldots \tau_1 (z) &\\
   &= \tau_i \ldots \tau_j \tau_{j+1} \tau_j \ldots 
    \tau_1 (z) &\\ 
   &= \tau_i \ldots \tau_{j+1} \tau_j \tau_{j+1} \ldots 
    \tau_1 (z) &\\ 
   &= \tau_i \ldots \tau_1 \tau_{j+1} (z) &\\
   &= \tau_i \ldots \tau_1 (1 \tensor \tau_j) (z)
   & \mbox{ for } j < i-1;\\
 \tau_{i-1(i)} \tau_{i-1} \ldots \tau_1 (z)
   &= \tau_i \tau_{i-1} \tau_i \tau_{i-1} \tau_{i-2} \ldots 
\tau_1 (z) &\\ 
   &= \tau_{i-1} \tau_i \tau_{i-1}^2 \tau_{i-2} \ldots 
\tau_1 (z) &\\ 
   &= \tau_{i-1} \tau_i \tau_{i-2} \ldots \tau_1 (z) &\\ 
   &= \tau_{i-1} \tau_{i-2} \ldots \tau_1 \tau_i (z) &\\
   &= \tau_{i-1} \ldots \tau_1 (1 \tensor \tau_{i-1}) (z);
   &\\
 \tau_{i(i)} \tau_i \ldots \tau_1 (z)
   &= \tau_i \tau_{i+1} \tau_i \ldots \tau_1 (z) &\\
   &= \tau_{i+1} \tau_i \tau_{i+1} \ldots \tau_1 (z) &\\
   &= \tau_{i+1} \tau_i \ldots \tau_1 \tau_{i+1}(z) &\\
   &= \tau_{i+1} \ldots \tau_1 (1 \tensor \tau_i) (z)
   &\\
 \end{array}$$
 where we used (\ref{equn2}) in the 3.~and 7.~equations.
 \end{pf}

 \begin{lma} \label{zetaimage}
  For all $z \in P^{n+1}(-1,\zeta)$ and all $f: P^2 \to P$ 
we have
 $$(P^{i-1} \tensor f \tensor P^{n-i}) \tau_{i-1} \ldots 
\tau_1(z) \in P^n(\zeta).$$ 
 \end{lma}
 
 \begin{pf}
 For all $\phi \in B_n$ and all $k = 1, \ldots, n$ we have
 $$\begin{array}{r @{\ } l}
 \tau_k^2 \phi &(P^{i-1} \tensor f \tensor P^{n-i}) 
\tau_{i-1} \ldots \tau_1(z) = \\ 
 &= \tau_k^2 (P^{j-1} \tensor f \tensor P^{n-j}) \phi_{(i)} 
\tau_{j-1} \ldots \tau_1(z) \\ 
 &= (P^{j-1} \tensor f \tensor P^{n-j}) \tau_{k(j)}^2 
\tau_{j-1} \ldots \tau_1 (1 \tensor \phi) (z) \\ 
 &= (P^{j-1} \tensor f \tensor P^{n-j}) \tau_{j-1} \ldots 
\tau_1 (1 \tensor \tau_k^2 \phi) (z) \\ 
 &= \phi (P^{i-1} \tensor f \tensor P^{n-i}) \tau_{i-1} 
\ldots \tau_1 (z) \\ 
 \end{array}$$
 hence $(P^{i-1} \tensor f \tensor P^{n-i}) \tau_{i-1} 
\ldots \tau_1(z) \in P^n(\zeta).$ 
 \end{pf}

 Now we can give the
 \begin{pf} of Proposition \ref{jacobidef1}:

 We first show that $P^{n+1}(\zeta) \subseteq P \tensor 
(P^n(\zeta)) \subseteq P^{n+1}$. Let $z = \sum_k z_{k,1} 
\tensor \ldots \tensor z_{k,n+1}$ be in $P^{n+1}(\zeta)$ 
with linearly independent $z_{k,1}$. Let $\phi, \tau_i \in 
B_n$ be given. Define $1 \tensor \phi \in B_{n+1}$ resp. $1 
\tensor \tau_i \in B_{n+1}$ by the operation of $\phi$ 
resp. $\tau_i$ on the factors $z_{k,2} \tensor \ldots 
\tensor z_{k,n+1}$, e.g. $1 \tensor \tau_i^{(n)} = 
\tau_{i+1}^{(n+1)}$. Then 
 $$\begin{array}{l @{\ } l}
 \sum z_{k,1} \tensor \phi^{-1} \tau_i^2 \phi (\sum z_{k,2} 
\tensor \ldots \tensor z_{k,n+1}) 
 &= (1 \tensor \phi^{-1} \tau_i^2 
\phi)(z) \\ 
 &= \sum_k z_{k,1} \tensor \zeta^2 \sum z_{k,2} \tensor 
\ldots \tensor z_{k,n+1}.\\ 
 \end{array}$$ 
 Since the $z_{k,1}$ are linearly independent, the terms 
$\sum z_{k,2} \tensor \ldots \tensor z_{k,n+1}$ are in 
$P^n(\zeta)$ hence $z \in P \tensor P^n(\zeta)$.

 Now we show that a factorization as given in the following 
diagram exists 
  $$\bfig 
 \putmorphism(0, 500)(1, 0)[P^{n+1}(\zeta)`P \tensor 
P^n(\zeta)`\iota]{800}1a 
 \putmorphism(0, 0)(1, 0)[P^2(-1)`P \tensor 
P.`\iota]{800}1a 
 \putmorphism(0, 500)(0, 1)[``1 \tensor \mbox{[.,.]} 
]{500}1r 
 \putmorphism(800, 500)(0, 1)[``1 \tensor 
\mbox{[.,.]}]{500}1r 
 \efig$$ 
 The morphism $1 \tensor [.,.]: P \tensor P^n(\zeta) \to P 
\tensor P$ is in $\C$. Consider the braiding $\tau: P 
\tensor P^n(\zeta) \to P^n(\zeta) \tensor P$. Since it is a 
natural transformation the diagram 
 $$\bfig 
 \putmorphism(0, 500)(1, 0)[P \tensor P^n(\zeta)`P^n(\zeta) 
\tensor P`\tau]{1200}1a 
 \putmorphism(0, 0)(1, 0)[P \tensor (P \tensor \ldots 
\tensor P)`(P \tensor \ldots \tensor P) \tensor 
P`\phi]{1200}1a 
 \putmorphism(0, 500)(0, 1)[``]{500}1r 
 \putmorphism(1200, 500)(0, 1)[``]{500}1r 
 \efig$$ 
 commutes with $\phi = \tau_n \ldots \tau_1 = 
\tau_{(P,P^n)}$, so $\tau(\sum_k z_{k,1} \tensor (\sum 
z_{k,2} \tensor \ldots \tensor z_{k,n+1})) = \tau_n \ldots 
\tau_1(\sum_k z_{k,1} \tensor z_{k,2} \tensor \ldots 
\tensor z_{k,n+1})$. Hence we get 
 \begin{equation}\label{equn3} 
 \tau(1 \tensor [.,.])(z) 
 = ([.,.] \tensor 1)\tau_n \ldots \tau_1(z) 
 \end{equation} 
 and similarly 
 \begin{equation}\label{equn4} 
 \tau([.,.] \tensor 1)(z) = (1 \tensor [.,.])\tau_1 \ldots 
\tau_n(z) 
 \end{equation} 
 for $z \in P^{n+1}(\zeta)$. This implies $\tau^2(1 \tensor 
[.,.])(z) = \tau([.,.] \tensor 1)\tau_n \ldots \tau_1(z) = 
(1 \tensor [.,.])\tau_1 \ldots \tau_n \tau_n \ldots 
\tau_1(z) = (1 \tensor [.,.]) \zeta^{2n}(z) = (-1)^2 (1 
\tensor [.,.])(z)$, so that $(1 \tensor [.,.])(z)$ is in 
$P^2(-1)$ and thus $\sum [z_{k,1},[z_{k,2}, \ldots, 
z_{k,n+1}]]$ is defined. 

 The second claim of the Proposition is proved in a 
symmetric way.
 \end{pf} 

 We continue with the
 \begin{pf} of Proposition \ref{jacobidef2}:

 We use (\ref{equn3}), (\ref{equn4}), and (\ref{equn1}) to 
get 
 $$\begin{array}{r @{\ } l} 
 \tau^2 (1 \tensor &[.,.]_n)(z) = \\ 
 &= (1 \tensor [.,.]_n)\tau_1 \ldots \tau_n \tau_n \ldots 
\tau_1 (z)\\ 
 &= (1 \tensor [.,.]_n)\tau_1 \ldots \tau_{n-1 } \tau_n^2 
\tau_{n-1}^{-1} \ldots \tau_1^{-1} \tau_1 \ldots \tau_{n-1} 
\tau_{n-1} \ldots \tau_1 (z) \\ 
 &= (1 \tensor [.,.]_n) \tau_n \ldots \tau_2 \tau_1^2 
\tau_2^{-1} \ldots \tau_n^{-1} \tau_1 \ldots \tau_{n-2} 
\tau_{n-1}^2 \tau_{n-2} \ldots \tau_1 (z) \\ 
 &= (1 \tensor [.,.]_n) (\tau_n \ldots \tau_2 \tau_1^2 
\tau_2^{-1} \ldots \tau_n^{-1}) (\tau_{n-1} \ldots \tau_2 
\tau_1^2 \tau_2^{-1} \ldots \tau_{n-1}^{-1}) \ldots 
(\tau_1^2) (z)\\ 
 &= (1 \tensor [.,.]_n) (z)\\ 
 \end{array}$$ 
 for all $z \in P^{n+1}(-1,\zeta)$ hence $(1 \tensor 
[.,.]_n)(z)$ is in $P^2(-1)$ and $[.,[.,.]_n]_2(z)$ is 
defined. 

 Now we prove that $[.,[.,.]_2,.]_n \tau_{i-1} \ldots 
\tau_1: P^{n+1}(-1,\zeta) \ni x \tensor y_1 \tensor \ldots 
\tensor y_n \mapsto [y_1, \ldots , [x, y_i], \ldots, y_n] 
\in P$ is well defined. Let $z \in P^{n+1}(-1,\zeta)$. 
Then we have 
 $\tau_1^2 \tau_2^{-1} \ldots \tau_i^{-1}(z) 
 = \tau_2^{-1} \ldots \tau_i^{-1}(z)$ 
 since $\tau_2^{-1} \ldots \tau_i^{-1} = 1 \tensor 
\tau_1^{-1} \ldots \tau_{i-1}^{-1} $. If we represent $y = 
\tau_2^{-1} \ldots \tau_i^{-1}(z) = \sum a_i \tensor b_i 
\in P^2 \tensor P^{n-1}$ in shortest form, then the set 
$\{b_i\}$ is linearly independent, so $\sum a_i \tensor b_i 
= \tau_1^2(\sum a_i \tensor b_i) = \sum \tau_1^2(a_i) 
\tensor b_i$, hence $\tau_1^2(a_i) = a_i$ and $y \in P^2(-
1) \tensor P^{n-1}$. So we get $\tau_{i-1} \ldots \tau_1 
\tau_i \ldots \tau_2(y) = \tau_{i-1} \ldots \tau_1(z) \in 
P^{i-1} \tensor P^2(-1) \tensor P^{n-i}$ and $(1 \tensor 
\ldots \tensor [.,.] \tensor \ldots \tensor 1)\tau_{i-1} 
\ldots \tau_1(z) \in P^n$ is defined. 

 By Lemma \ref{zetaimage} we have $(1 \tensor \ldots 
\tensor [.,.]_2 \tensor \ldots \tensor 1)(z) \in 
P^n(\zeta)$, so that\break $[.,[.,.]_2,.]_n \tau_{i-1} \ldots 
\tau_1(z)$ is well defined.
 \end{pf} 


\end{document}